\documentclass{article}[12pt,notitlepage]
\usepackage{amsmath}
\usepackage{a4wide}
\usepackage{amsfonts}
\usepackage[pdftex]{graphicx}
\usepackage{float}
\usepackage[super, sort]{natbib}
\bibpunct{(}{)}{;}{a}{,}{,} 

\newtheorem{theorem}{Theorem}
\newtheorem{lemma}{Lemma}
\newtheorem{definition}{Definition}
\newtheorem{condition}{Condition}
\numberwithin{theorem}{section}
\numberwithin{lemma}{section}
\numberwithin{definition}{section}
\numberwithin{condition}{section}
\numberwithin{figure}{section}

\begin{document}

\title{Non-invertibility Criteria in Some Heteroscedastic Models}



\author{Alexey Sorokin\footnote{135 Hermitage Waterside, Thomas More Street, London, United Kingdom. Tel: +44 785 242 8329. Email address: alexey.a.sorokin@gmail.com.} \\ Knight Capital Europe Limited}
\date{}
\maketitle

\begin{abstract}
In order to calculate the unobserved volatility in conditional heteroscedastic
time series models, the natural recursive approximation is very often used.
Following \cite{StraumannMikosch2006}, we will call the
model \emph{invertible} if this approximation (based on true parameter vector)
converges to the real volatility. Our main results are necessary
and sufficient conditions for invertibility.
We will show that the stationary
GARCH($p$, $q$) model is always invertible, but certain types
of models, such as EGARCH of \cite{Nelson1991} and
VGARCH of \cite{EngleNg1993} may indeed be non-invertible.
Moreover, we will demonstrate it's
possible for the pair (true volatility, approximation)
to have a non-degenerate stationary distribution.
In such cases, the volatility
estimate given by the recursive approximation with the true
parameter vector is inconsistent.
\end{abstract}


\textit{Journal Of Economic Literature Clasification:} C14, C58.

\section{Introduction}

Since their introduction in the seminal papers of \cite{Engle1982}
and \cite{Bollerslev1986}, there has been a remarkable amount
of interest in heteroscedastic models coming from researchers in
econometrics and statistics and from financial market
practitioners. Since then, many other similar
models have been introduced to more accurately reflect different qualities of real data.
For example, various authors have proposed a vast number of volatility
expressions with asymmetries or thresholds, non-stationarities such as
unit roots or time-varying parameters, or have considered multidimensional
extensions and non-regular time intervals. We refer to two of more
recent reviews \citet{AndersenEtAl2009} and 
\cite{FrancqZakoian2010} which discuss these
and other examples and provide the necessary references.

Heteroscedastic models are used to analyze and forecast the
unobserved volatility of
different financial time series. It is common to
construct predictions for current and future volatility values using a
natural recursive volatility approximation. Let us recall its
definition for a general conditional heteroscedastic model.
We define such model as a solution of
\begin{equation}\label{equ_garch_general}
\begin{split}
y_t = \sigma_t \varepsilon_t,\\
\sigma^2_t = H(\theta, y_{t-1}, \ldots, y_{t-p},
\sigma^2_{t-1}, \ldots \sigma^2_{t-q}), \quad t \in \mathbb{Z},
\end{split}
\end{equation}
where $\{y_t\}$ is the return process,
 $\{\sigma_t\}$ is the volatility process, $H$ is a known function,
$p, q \in \mathbb{N} \cup \{0\}$,
$\theta$ is the parameter vector taken from a set $\Theta \subseteq \mathbb{R}^l$
and $\{\varepsilon_t\}$ are i.i.d. random variables on a probability 
space $(\Omega, \mathcal{F}, P)$. More general
versions of (\ref{equ_garch_general}) are of course possible. Assume that
$y_{1-p}, \ldots, y_n$ are observed, fix $s^2 \in \mathbb{R}^+$ and define
the recursive volatility approximation by
\begin{equation}\nonumber
\begin{split}
\hat{\sigma}^2_t(\theta') = s^2, \quad t = 1 - q, \ldots, 0,\\
\hat{\sigma}^2_t(\theta') = H(\theta', y_{t-1}, \ldots, y_{t-p},
\hat{\sigma}^2_{t-1}(\theta'), \ldots \hat{\sigma}^2_{t-q}(\theta')), \quad t > 0.
\end{split}
\end{equation}
for any $\theta' \in \Theta$. The quantity
\begin{equation}\label{equ_vol_forecast}
\hat{\sigma}^2_n(\hat{\theta}_n)
\end{equation}
is often used as an estimate of the volatility value at $n$, where $\hat{\theta}_n$
is any consistent estimator of $\theta$. There are other possible 
definitions of $\hat{\sigma}^2_t$. For example, \cite{BaillieEtAl1996} suggest using the sample mean of $y_{1-p}^2, \ldots, y_n^2$ as a
starting point instead of $s^2$.

In order for approximations (\ref{equ_vol_forecast}) to be consistent, it
is natural to at least expect that the approximation
$\hat{\sigma}^2_{t}(\theta)$ (using the \emph{true} parameter vector)
converges to the unobserved value $\sigma^2_t$:
\begin{equation}\label{equ_vol_converges}
|\hat{\sigma}^2_{t}(\theta) - \sigma^2_t| \stackrel{P} \rightarrow 0, \quad t \to \infty.
\end{equation}
To describe models with property (\ref{equ_vol_converges}), we'll re-use
the following definition from \cite{StraumannMikosch2006}:
\begin{definition}\label{def_invertibility}
Assume for all $\theta \in \Theta$ the model (\ref{equ_garch_general})
has a strictly stationary
solution and (\ref{equ_vol_converges}) is true for any initial point $s^2$.
Then we say the model has the property of \emph{invertibility}.
\end{definition}

Note that many estimator types
and tests are very commonly considered under the assumption of invertibility.
Indeed, such estimators as the popular
quasi-maximum likelihood
estimator (QMLE) of \cite{LeeHansen1994} and
more recent
GMM-type robust estimators of \cite{Boldin2000} or
minimum distance estimators as in \cite{Sorokin2004}
as well as many others are based on 
residuals $\{y_t/\sigma_t(\theta')\}$. Their asymptotic
properties are commonly established assuming convergence of residuals,
therefore under invertibility. The same assumption is used for popular
tests types such as goodness-of-fit and dimensionality tests
(see \cite{FrancqZakoian2010})
and structural break tests (\cite{Boldin2002},
\cite{HorvathTeyssiere2001}). Therefore, the issue of non-invertibility is
key for such procedures. The notable exception
are estimators based on autocovariance function
such as Whittle estimator applied in \cite{GiraitisRobinson2000} and \cite{Zaffaroni2009}.

For linear time series models such as ARMA the notion of invertibility is classic, see e.g. \cite{BrockwellDavis1991}. For non-linear models, however,
there seem to be more than one possible definition.
The definition \ref{def_invertibility} is based on that of \cite{GrangerAndersen1978}.
For detailed discussion and some further references, see \cite{StraumannMikosch2006}, subsection 3.2. 
A related definition of invertibility, in particular for bilinear models, was
considered previously in \cite{Tong1990}. According to Tong's
definition, the model (\ref{equ_garch_general}) is called invertible
if $\varepsilon_t$ is a.s. a measurable function of $(\ldots, y_{t-2}, y_{t-1}, y_t)$
for any $t \in \mathbb{Z}$. We'll refer to this notion as global 
invertibility (this is in contrast to local invertibility, see the next paragraph). 
Recently, there have been some further results on invertibility of non-linear models. 
The threshold MA models were investigated in \cite{LingTong2005} 
and \cite{LingEtAl2007},
and non-linear ARMA (NLARMA) models were discussed in \cite{ChanTong2009}.

The main goal of this paper is to show that some heteroscedastic models
(for example, EGARCH of \cite{Nelson1991} and VGARCH of
\cite{EngleNg1993}) may in general be non-invertible,
even if the model has a strictly stationary solution. Moreover, we will 
provide necessary and sufficient conditions for invertibility for such models. 
We will demonstrate that if a model is 
locally invertible (see Definition \ref{def_proper_unidentifiability} below, also 
\cite{ChanTong2009} for the similar definition), it is also invertible. On the other
hand, we'll prove that under a slightly more restricting condition than a
lack of local invertibility the model is also non-invertible.
Moreover, we will show that on an extended probability space 
there exists a random
variable $\xi$, independent with $\sigma\{\varepsilon_t, \; t \geq 0\}$,
such that for $s^2 = \xi$ the process
$$
(\sigma^2_t, \hat{\sigma}^2_t(\theta)), \quad t \geq 0
$$
is stationary but 
$$
P(\sigma^2_t = \hat{\sigma}^2_t(\theta)) = 0, \quad t \geq 0.
$$
To the best of our knowledge, this is the first use of the notion of local invertibility in
the literature on heteroscedastic models. 
Its usefulnes is that it's much easier to check that invertibility itself, this will be discussed below.
The used method is general and may be applied to many other heteroscedastic models. 


The rest of the paper is organized as follows: in the subsection
\ref{sec_stat_inv} we provide necessary definitions and discuss
non-invertibility in more detail. 
The main results of the paper are given
 in the section \ref{sec_main_results}. 
As a fact of independent interest and
an illustration of non-invertibility, we provide the results of numerical
simulations showing that it is possible to have deterministic chaos
(in the sense of \cite{Devaney1989}) in EGARCH in
subsection \ref{sec_egarch_chaos}.
Finally, section
\ref{sec_proofs} contains proofs.


\section{Prelimiaries}\label{sec_preliminaries}

\subsection{Stationarity and Invertibility}\label{sec_stat_inv}

Let us recall the definition of the GARCH($p$, $q$) model of
\cite{Bollerslev1986}. It is defined as a solution of particular
case of (\ref{equ_garch_general}) with $l = p + q + 1$ and
\begin{equation}\label{equ_garch}
\begin{split}
y_t = \sigma_t \varepsilon_t,\\
\sigma_t^2 = \alpha_0 + \alpha_1 y^2_{t-1} +
\ldots + \alpha_p y^2_{t-p} + \beta_1 \sigma^2_{t-1} + \ldots + \beta_q \sigma^2_{t-q},
\quad t \in \mathbb{Z},
\end{split}
\end{equation}
where
$$
\Theta = \{\alpha_0 > 0, \quad \alpha_i \geq 0, \; i = 1, \ldots , p,
\quad \beta_j \geq 0, \; j = 1, \ldots , q\}.
$$
We'll be interested in the strictly stationary solution of (\ref{equ_garch}).
The necessary and sufficient condition for the existence and uniqueness
of such a solution was found in \cite{BougerolPicard1992b}.
For simplicity, let's consider a case $p = q = 1$. Then the condition reads
\begin{condition}\label{cond_stat_garch}
$$
E \log(\alpha_1 \varepsilon_0^2 + \beta_1) < 0.
$$
\end{condition}
Here and in the rest of the paper we denote the natural logarithm $\log$.

It is easy to see that in GARCH($1$, $1$) model strict
stationarity implies invertibility. Indeed,
Condition \ref{cond_stat_garch} implies that $\beta_1 < 1$. Then
$$
\hat{\sigma}^2_t(\theta) - \sigma^2_t =  \alpha_0 + \alpha_1 y^2_{t-1} +
\beta_1 \hat{\sigma}^2_{t-1}(\theta) - (\alpha_0 + \alpha_1 y^2_{t-1} +
\beta_1 \sigma^2_{t-1}(\theta)) = \beta_1 (\hat{\sigma}^2_{t-1}(\theta) - \sigma^2_{t-1}),
$$
therefore
$$
|\hat{\sigma}^2_t(\theta) - \sigma^2_t| \stackrel{a.s.} \rightarrow 0,
 \quad t \to \infty,
$$
implying invertibility. Using the multidimensional version of Condition
\ref{cond_stat_garch}, one can check that the same implication is also
true for the general GARCH($p$, $q$).

However, the stationarity does not imply invertibility for other
heteroscedastic models. Notably, it turns out this relation fails
for the EGARCH model of \cite{Nelson1991} and the VGARCH model of \cite{EngleNg1993}.
For simplicity, we'll consider only 1-dimensional models of those types.

The general model (\ref{equ_garch_general})
for the 1-dimensional case may be written as
\begin{equation}\label{equ_garch_general_1dim}
\begin{split}
y_t = \sigma_t \varepsilon_t,\\
\sigma^2_t = H(\theta, y_{t-1}, \sigma^2_{t-1}),
\quad t \in \mathbb{Z}.
\end{split}
\end{equation}
\begin{definition}\label{def_egarch}
EGARCH is defined as a solution of
\begin{equation}\label{equ_egarch}
\begin{split}
y_t = \sigma_t \varepsilon_t,\\
\log \sigma_t^2 = \alpha + \gamma |\varepsilon_{t-1}| + 
\delta \varepsilon_{t-1} + \beta \log \sigma^2_{t-1}, \quad t \in \mathbb{Z},
\end{split}
\end{equation}
where $\theta = (\alpha, \beta, \gamma, \delta)$ is the parameter vector and $\{\varepsilon_t\}$ are i.i.d.
\end{definition}

\begin{definition}\label{def_vgarch}
VGARCH is defined as a solution of
\begin{equation}\label{equ_vgarch}
\begin{split}
y_t = \sigma_t \varepsilon_t,\\
\sigma_t^2 = \alpha + \gamma (\varepsilon_{t-1} - \delta)^2 + \beta \sigma^2_{t-1},
\quad t \in \mathbb{Z},
\end{split}
\end{equation}
where $\theta = (\alpha, \beta, \gamma, \delta)$ is the parameter vector and $\{\varepsilon_t\}$ are i.i.d.
\end{definition}
In order to make the models meaningful, we need to enforce certain restriction on parameter values.
For EGARCH, we assume that $\gamma \geq |\delta|$ (see \cite{StraumannMikosch2006}, 
p. 20 for discussion of why this is reasonable, note somewhat different notation) and that $\beta > 0$. 
For VGARCH, we demand that $\alpha > 0, \; \beta > 0, \gamma \geq 0$.

It follows from the Theorem 2.5 from \cite{BougerolPicard1992a} that a necessary and sufficient conditions for the
existence of a strictly stationary solution to (\ref{equ_egarch})
and (\ref{equ_vgarch}) are given by
\begin{condition}\label{cond_stat_egarch}
$$
E \log \left |\alpha + \gamma |\varepsilon_0| + \delta \varepsilon_0 \right| < \infty, \quad \beta < 1.
$$
\end{condition}
\begin{condition}\label{cond_stat_vgarch}
$$
E \log \left (\alpha + \gamma (\varepsilon_0 - \delta)^2 \right ) < \infty, \quad \beta < 1.
$$
\end{condition}
respectively. We will assume these conditions to hold throughout the 
rest of the paper, and denote the stationary solution of either (\ref{equ_egarch})
or (\ref{equ_vgarch}) by $\{\overline{\sigma}_t\}$.

The issue of invertibility is, however, more complicated.
First, let's consider the explicit expressions for $\hat{\sigma}^2(\theta)$.
Since the sequence of random variables
$\{\varepsilon_t\}$ is unobservable, it needs to be approximated.
Commonly this is done using the residuals
$$
\hat{\varepsilon}_t(\theta') = y_t / \hat{\sigma}_t(\theta')
$$
For EGARCH, we define $\hat{\sigma}^2(\theta')$ recursively as a solution of
\begin{equation}\label{equ_egarch_sigmahat}
\begin{split}
\log \hat{\sigma}^2_t(\theta') = s^2, \quad t = 0,\\
\log \hat{\sigma}^2_t(\theta') = \alpha' + 
\beta' \log \hat{\sigma}^2_{t-1}(\theta') + \gamma'
|\hat{\varepsilon}_{t-1} (\theta')| + \delta' \hat{\varepsilon}_{t-1} (\theta'), \quad t > 0.
\end{split}
\end{equation}
For VGARCH $\hat{\sigma}^2(\theta')$ are given by
\begin{equation}\label{equ_vgarch_sigmahat}
\begin{split}
\log \hat{\sigma}^2_t(\theta') = s^2, \quad t = 0,\\
\hat{\sigma}^2_t(\theta') = \alpha' + 
\beta' \hat{\sigma}^2_{t-1}(\theta') + \gamma'
(\hat{\varepsilon}_{t-1}(\theta') - \delta')^2, \quad t > 0.
\end{split}
\end{equation}
In both cases, we have
$$
\hat{\sigma}^2_t(\theta') = f_{t-1}(\hat{\sigma}^2_{t-1}(\theta')), \quad t > 0,
$$
where $\{f_{t}(\cdot)\}$ is a sequence of random transformations $\mathbb{R}^+ \to \mathbb{R}^+$, defined for the general model (\ref{equ_garch_general_1dim}) by
$$
f_{t}(x) = H(\theta, y_{t}, x), \quad t \in \mathbb{Z}, \; x \in \mathbb{R}^+.
$$
Crucially, the random transformations defining $\{\sigma^2_t\}$
(cf. (\ref{equ_egarch}), (\ref{equ_vgarch})) and those defining
$\{\hat{\sigma}^2_t(\theta)\}$ (cf. (\ref{equ_egarch_sigmahat}), (\ref{equ_vgarch_sigmahat}))
are different. Indeed, for the general model
(\ref{equ_garch_general_1dim}) we have
\begin{equation}\label{equ_sigma_sigmahat_system}
\begin{cases}
\sigma^2_t = H(\theta, \sigma_{t-1} \varepsilon_{t-1}, \sigma^2_{t-1}),\\
\hat{\sigma}^2_t(\theta) = H(\theta, \sigma_{t-1} \varepsilon_{t-1},
\hat{\sigma}^2_{t-1}(\theta))
\end{cases}
\quad t \geq 1.
\end{equation}
Note that the right-hand side of the second equation in
(\ref{equ_sigma_sigmahat_system}) depends not only
on $\hat{\sigma}^2_{t-1}(\theta)$, but also on $\sigma^2_{t-1}$.

For a general heteroscedastic model invertibility
holds if the top Lyapunov exponent
of the corresponding sequence $\{f_{t}(\cdot)\}$ is
negative (Proposition 3.7 of \cite{StraumannMikosch2006}).
Let us recall the definition of Lyapunov exponent of a family of
transformations:
\begin{definition}\label{def_lyapunov_exp}
Assume $\{f_t(\cdot)\}_{t \in \mathbb{Z}}$ is a stationary ergodic
family of $\mathbb{R} \to \mathbb{R}$ functions, such that
$$
E \log \Lambda(f_t) < \infty,
$$
where $\Lambda$ is the Lipchitz norm,
$$
\Lambda(f) := \sup_{x \neq y} \frac{|f(x) - f(y)|}{|x - y|}.
$$
Then the quantity
$$
\lambda = \inf \left \{E \left ( \frac{1}{t + 1}
\log \Lambda(f_t \circ \ldots \circ f_0) \right), t \in \mathbb{N} \right \}
$$
is called the top Lyapunov exponent.
\end{definition}
However, estimating top Lyapunov exponent for nonlinear systems
is not straightforward (and therefore checking whether the sufficient 
condition for invertibility holds is hard, p. 21 of \cite{StraumannMikosch2006}). 
Note that the second equation of (\ref{equ_sigma_sigmahat_system})
could be nonlinear in $\hat{\sigma}^2_{t-1}(\theta)$ even if the first one is linear!
Therefore, the results of \cite{BougerolPicard1992a} are not applicable.
To the best of our knowledge, necessary and sufficient conditions for
invertibility of such nonlinear models are not known.

We will consider a simple lower
bound for the Lyapunov exponent which will be useful for the construction of a criteria for invertibility.
From now on, we will write $\hat{\sigma}^2_t$
instead of $\hat{\sigma}^2_t(\theta)$ for brevity.
\begin{definition}\label{def_stability_coeff}
Assume that the equation (\ref{equ_garch_general_1dim}) has a unique stationary
solution $\{\overline{\sigma}^2_t\}_{t \in \mathbb{Z}}$. Assume also that
for any fixed $\theta \in \Theta$ the function $H(\theta, \cdot, \cdot)$ 
is continuously differentiable. Denote
$$
\Lambda^y_t = \left | \left . \frac{\partial}{\partial \hat{\sigma}^2_0}
\right |_{\hat{\sigma}^2_0 = \overline{\sigma}^2_0} (f_t \circ \ldots \circ f_0) \right |,
$$
$$
\lambda^y = \inf \left \{E \left ( \frac{1}{t + 1}
\log \Lambda^y_t \right), t \in \mathbb{N} \right \}.
$$
We will refer to $\lambda^y$ as \emph{stability coefficient} and call a
heteroscedastic model locally invertible if $\lambda^y < 0$ and
locally non-invertible otherwise.
\end{definition}
A similar definition for non-linear ARMA models was considered in 
\cite{ChanTong2009}.

Unlike top Lyapunov exponent, $\lambda^y$ may be written down explicitly
using the chain differentiation rule. Denote
$$
H'_x(\theta, y, x) = \frac{\partial H(\theta, y, x) }{\partial x},
$$
then
$$
\Lambda^y_t = \left | \left . \frac{\partial}{\partial \hat{\sigma}^2_0}
\right |_{\hat{\sigma}^2_0 = \overline{\sigma}^2_0} (f_t \circ \ldots \circ f_0) \right | =
\prod_{k = 0}^t \left | H'_x(\theta, \overline{\sigma}_k \varepsilon_k, \overline{\sigma}^2_k)\right |,
$$
therefore
$$
\log \Lambda^y_t = \sum_{k = 0}^t \log \left | H'_x(\theta, \overline{\sigma}_k
 \varepsilon_k, \overline{\sigma}^2_k) \right |.
$$
If the process $\{\overline{\sigma}^2_t\}$ is ergodic and
$$
E \log \left | H'_x(\theta, \overline{\sigma}_0 \varepsilon_0,
\overline{\sigma}^2_0) \right | < \infty,
$$
due to Birkhoff-Khinchin ergodic theorem
\begin{equation}\label{equ_lambda_y_formula}
\lambda^y = E \log \left | H'_x(\theta, \overline{\sigma}_0 \varepsilon_0,
\overline{\sigma}^2_0) \right |.
\end{equation}
Also note that since
$$
\Lambda^y_t \leq \Lambda(f_t \circ \ldots \circ f_0),
$$
it must also be that
$$
\lambda^y \leq \lambda.
$$
The quantity $\lambda^y$ has an intuitive interpretation. It indicates
whether the behaviour of $\hat{\sigma}^2_t$ is convergent or explosive in
the neighborhood of the true volatility value $\sigma^2_t$ (this will be formalized in Theorems 
\ref{thm_stat_distr_egarch} - \ref{thm_lln_noninv_vgarch}).
For the invertibility of (\ref{equ_garch_general_1dim}),
three cases are possible:

1. If $\lambda^y \leq \lambda < 0$, $\hat{\sigma}^2_t$ is stable and
invertibility takes place (see \cite{StraumannMikosch2006}, Theorem 2.8).

2. If $\lambda^y \leq 0 \leq \lambda$, Theorem 2.8 of \cite{StraumannMikosch2006} is
not valid anymore. The behaviour of $\hat{\sigma}^2_t$ is "stable" in
the neighborhood of the true volatility value $\sigma^2_t$.

3. If $0 < \lambda^y \leq \lambda$, the behaviour of $\hat{\sigma}^2_t$ is explosive in
the neighborhood of the true volatility value $\sigma^2_t$. 


\begin{definition}\label{def_proper_unidentifiability}
Assume that the system  (\ref{equ_sigma_sigmahat_system})
has a stationary solution $(\overline{\sigma}^2_t, \hat{\overline{\sigma}}^2_t)_{t \geq 0}$
on an extended probability space such that
$$
P(\overline{\sigma}^2_0 = \hat{\overline{\sigma}}^2_0) = 0, \quad \geq 0.
$$
Then we call the model (\ref{equ_garch_general_1dim})
\emph{properly non-invertible}.
\end{definition}

It is interesting that the unstable assumption of case 3 is
(under mild technical conditions) sufficient for the
proper non-invertibility for EGARCH and VGARCH models. 
On the other hand, in case 2 if we only assume that $\lambda^y < 0$, 
both models are invertible. 
These are the main
results of this paper, to be considered in the following section.
In addition, we show that in both cases for any starting point $s^2$ the
sample distribution of
$$
(\sigma^2_t, \hat{\sigma}^2_t), \quad t = 0, \ldots, n
$$
converges to the stationary one. Together with proper
non-invertibility it implies that the
approximation $\hat{\sigma}^2_t$ is inconsistent if $\lambda^y > 0$.


\section{Main Results}\label{sec_main_results}

Before we formulate the main results, let's introduce the conditions
to be used.
First one is a regularity condition.
\begin{condition}\label{cond_regul}
The distribution of $\varepsilon_0$ is absolutely continuous with respect
to Lebesgue measure on $\mathbb{R}$ with continuous density $g(x)$,
$E \varepsilon_0 = 0$, $E \varepsilon_0^2 = 1$ and
$$
g(x) > 0 \; \forall x \in \mathbb{R}, \qquad \sup_{x \in \mathbb{R}} (|x| + 1) g(x) < \infty.
$$
\end{condition}
Condition \ref{cond_regul} may probably be somewhat relaxed, although a 
condition similar to $g(x) > 0$ is essential in establishing the $\psi$-irreducibility of 
relevant Markov chains, see Proofs section.


We will formulate the assumption on $\lambda^y$ using the equation
(\ref{equ_lambda_y_formula}). Denote $\lambda^y_E$ and $\lambda^y_V$
stability coefficients for EGARCH and VGARCH respectively. 
Note that for EGARCH due to (\ref{equ_egarch_sigmahat})
$$
H^E(\theta, y, x) = \exp \left \{\alpha + \beta \log x + 
(\gamma |y| + \delta y ) x^{- 1/2} \right \}, \quad \left ( H^E \right)'_x(\theta, y, x) =
 H^E(\theta, y, x) \left ( \beta x^{-1} - (\gamma |y| + \delta y ) x^{-3/2} / 2 \right ),
$$
and by definition of the stationary solution $\overline{\sigma}_t$
$$ 
\lambda^y_E = E \log \left |H^E(\theta,  \overline{\sigma}_0 \varepsilon_0, 
\overline{\sigma}^2_0) \right | - E \log |\overline{\sigma}^2_0| + 
E \log \left |\beta - (\gamma |\overline{\sigma}_0 \varepsilon_0| + 
\delta \overline{\sigma}_0 \varepsilon_0 ) \overline{\sigma}^{-1}_0  / 2 \right | = 
$$
$$
E \log \left |\beta - 2^{-1} (\gamma |\varepsilon_0| + \delta \varepsilon_0 ) \right |.
$$
For VGARCH due to (\ref{equ_vgarch_sigmahat})
$$
H^V(\theta, y, x) = \alpha + \beta x + \gamma (y x^{-1/2} - \delta)^2, \quad
\left ( H^V \right )'_x(\theta, y, x) = \beta - \gamma y x^{-3/2} (y x^{-1/2} - \delta)
$$
and
$$ 
\lambda^y_V = E \log \left | \beta - \gamma \overline{\sigma}_0 \varepsilon_0
\overline{\sigma}_0^{-3} (\overline{\sigma}_0 \varepsilon_0 \overline{\sigma}_0^{-1} -
 \delta) \right | = E \log \left | \beta - \gamma \varepsilon_0 \overline{\sigma}_0^{-2} (\varepsilon_0 - \delta) \right |.
$$

\subsection{Invertibility criteria}
Now we will formulate the main results:
\begin{theorem}[Proper non-invertibility for EGARCH]\label{thm_stat_distr_egarch}
Assume the Conditions \ref{cond_stat_egarch} and \ref{cond_regul} hold and $\lambda^y_E > 0$. Then the system (\ref{equ_sigma_sigmahat_system})
with $H = H^E$ has a stationary solution
$\{(\overline{\sigma}^2_t, \hat{\overline{\sigma}}^2_t)\}_{t \geq 0}$ on an 
extended probability space such that
$$
P(\overline{\sigma}^2_0 = \hat{\overline{\sigma}}^2_0) = 0.
$$

\end{theorem}
\begin{theorem}[Proper non-invertibility for VGARCH]\label{thm_stat_distr_vgarch}
Assume the Conditions \ref{cond_stat_vgarch} and \ref{cond_regul} hold and $\lambda^y_V > 0$. Then the system (\ref{equ_sigma_sigmahat_system})
with $H = H^V$ has a stationary solution
$\{(\overline{\sigma}^2_t, \hat{\overline{\sigma}}^2_t)\}_{t \geq 0}$ on an 
extended probability space such that
$$
P(\overline{\sigma}^2_0 = \hat{\overline{\sigma}}^2_0) = 0.
$$
\end{theorem}
Besides the proper non-invertibility, another important question is
the joint behaviour of the true volatility $\sigma^2_t$ and
its approximation $\hat{\sigma}^2_t$. As it turns out,
the sample distribution converges to the stationary one:
\begin{theorem}[LLN for $\{(\sigma^2_t, \hat{\sigma}^2_t)\}$ in EGARCH]
\label{thm_lln_noninv_egarch}
Assume the conditions \ref{cond_stat_egarch} and \ref{cond_regul} are
fulfilled. 

i) If $\lambda^y_E > 0$, the process $\{(\overline{\sigma}^2_t, \hat{\overline{\sigma}}^2_t)\}$,
the stationary solution of (\ref{equ_sigma_sigmahat_system})
with $H = H_E$, is ergodic and for any $B \in \mathcal{B}(\mathbb{R}^2)$ and starting point $s^2$
$$
n^{-1} \sum_{t = 0}^{n - 1} I \{(\sigma^2_t, \hat{\sigma}^2_t) \in B\}
\stackrel{a.s.} \rightarrow
P ((\overline{\sigma}^2_0, \hat{\overline{\sigma}}^2_0) \in B), \quad n \to \infty,
$$
where $I\{\cdot\}$ is the indicator function.

ii) If $\lambda^y_E < 0$, for any starting point $s^2$ and $\mu > 0$
$$
P (|\sigma^2_0 - \hat{\sigma}^2_0| \geq \mu) \to 0.
$$
\end{theorem}
\begin{theorem}[LLN for $\{(\sigma^2_t, \hat{\sigma}^2_t)\}$ in VGARCH]
\label{thm_lln_noninv_vgarch}
Assume the conditions \ref{cond_stat_egarch} and \ref{cond_regul} are
fulfilled. 

i) If $\lambda^y_V > 0$, the process
 $\{(\overline{\sigma}^2_t, \hat{\overline{\sigma}}^2_t)\}$,
the stationary solution of (\ref{equ_sigma_sigmahat_system})
with $H = H_V$,
is ergodic and for any $B \in \mathcal{B}(\mathbb{R}^2)$ and starting point $s^2$
$$
n^{-1} \sum_{t = 0}^{n - 1} I \{(\sigma^2_t, \hat{\sigma}^2_t) \in B\}
\stackrel{a.s.} \rightarrow
P ((\overline{\sigma}^2_0, \hat{\overline{\sigma}}^2_0) \in B), \quad n \to \infty.
$$
ii) If $\lambda^y_V < 0$, for any starting point $s^2$ and $\mu > 0$
$$
P (|\sigma^2_0 - \hat{\sigma}^2_0| \geq \mu) \to 0.
$$
\end{theorem}
Theorems \ref{thm_stat_distr_egarch} - \ref{thm_lln_noninv_vgarch} indicate
that for some parameter vectors finite sample volatility approximations
do not converge to the actual volatility as the sample size increases.
In particular, the standard recursive volatility estimator is inconsistent (this is easily checked by setting 
$$
B = \{(\sigma^2, \hat{\sigma}^2) \in \mathbb{R}^2_+: |\sigma^2 - \hat{\sigma}^2| \geq \mu\}
$$
for a $\mu > 0$ and applying Theorems \ref{thm_lln_noninv_egarch} and \ref{thm_lln_noninv_vgarch}).

The general heteroscedastic model (\ref{equ_garch_general_1dim}) considered
in this paper can be re-written as a particular case of NLARMA
model (see \cite{ChanTong2009}, Equation (3)). The issue of invertibility
then also becomes
a special case of invertibility of NLARMA. However, the latter paper
discusses only local invertibility for such models. Therefore, our results
complement \cite{ChanTong2009} since we
establish that local invertibility is (except the case $\lambda^y = 0$ )
equivalent to invertibility for particular cases of
NLARMA. It would be very interesting to investigate global invertibility 
of heteroscedastic models such as EGARCH and VGARCH. It easy to see that
a stationary invertible model is necessary globally invertible. Indeed,
denote for any $t \in \mathbb{Z}$, $k \geq 1$
$$
\hat{\sigma}^2_{t, k, s^2} = H(\theta, y_{t-1}, 
H(\theta, y_{t-2}, \ldots H(\theta, y_{t-k}, s^2) \ldots ),
$$
so that $\hat{\sigma}^2_{t} = \hat{\sigma}^2_{t, t, s^2}$. Then due to
invertibility 
$$
\hat{\sigma}^2_{t, k, s^2} \stackrel{P} \to \sigma_t^2
$$
as $k \to \infty$ and for some sub-sequence $\{k_n\}$ 
$$
\hat{\sigma}^2_{t, k_n, s^2} \stackrel{a.s.} \to \sigma_t^2
$$
as $n \to \infty$. The reverse relation seems more complicated and is 
left for future research. 

Another problem (which also will not be considered here)
is the behaviour of tests and estimates unver non-invertibility. Commonly, 
consistency and asymptotic normalilty are established under invertibility 
assumption. However, it was shown in \cite{Zaffaroni2009} that the Whittle
estimator (which is not residual-based)  
remains asymptotically normal under non-invertibility too. It will be interesting 
to see what happens to common estimators such as QMLE.

\begin{figure}[h!]
\begin{center}
\includegraphics[scale=0.34]{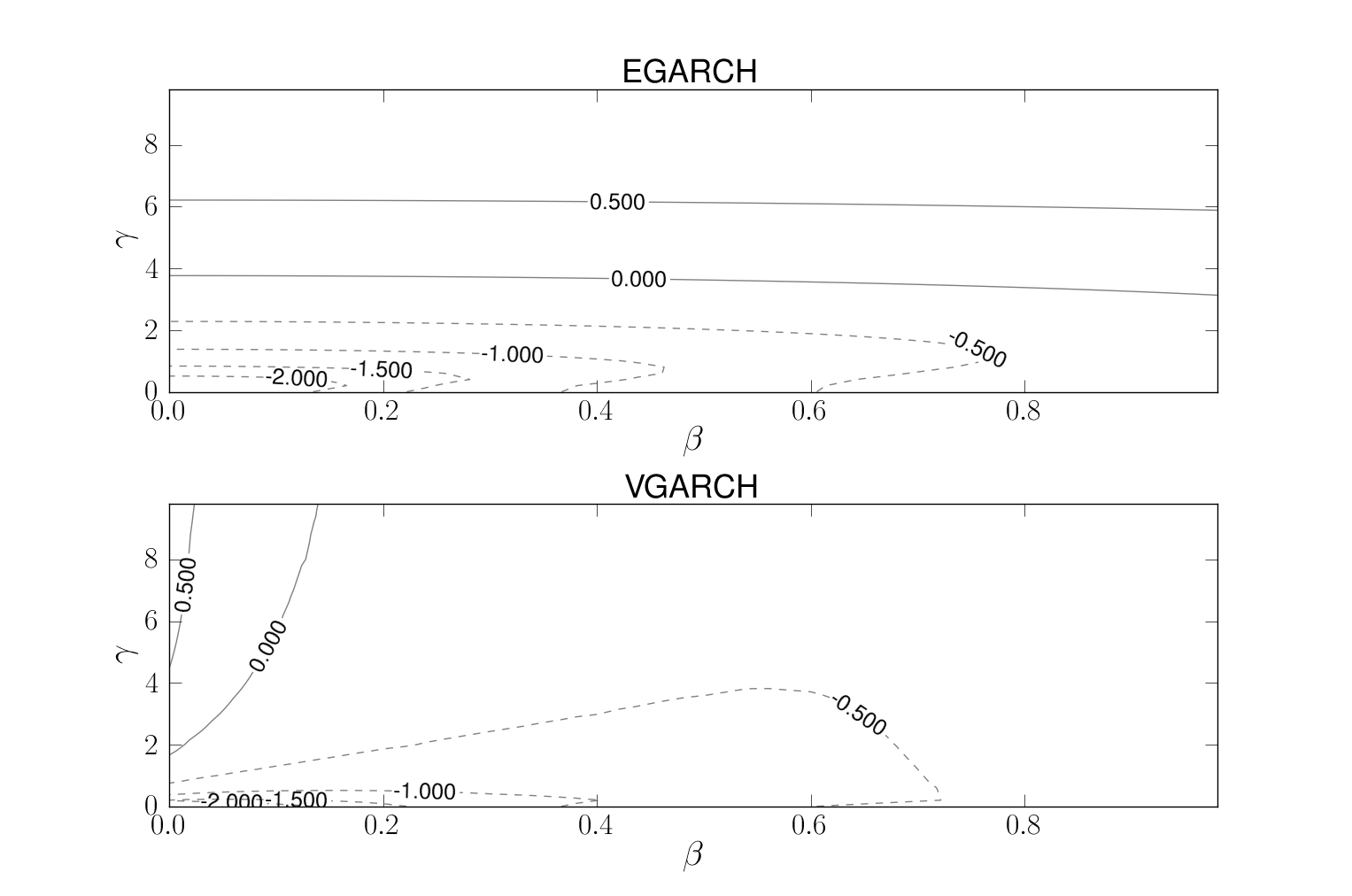}
\end{center}
\caption{Heatmap of $\lambda^y_E$ and $\lambda^y_V$. Positive values correspond to non-invertible models. For EGARCH, we used $\alpha = 0.1$, 
$\beta = 0.25$, $\gamma = 5.4$ and standard normal $\varepsilon_1$.
For VGARCH, we used $\alpha = 0.001$, 
$\beta = 0.01$, $\delta = -0.3$, $\gamma = 1$ and $\varepsilon_1$
was an appropriately scaled mixture of $\xi_1$ 
and $-\xi_2$, where both $\xi_1$ and $\xi_2$ had exponential distributions
 with parameters $1$ and $4$.}
\label{fig_lambda_y}
\end{figure}

\begin{figure}[H]

\begin{center}
\includegraphics[scale=0.34]{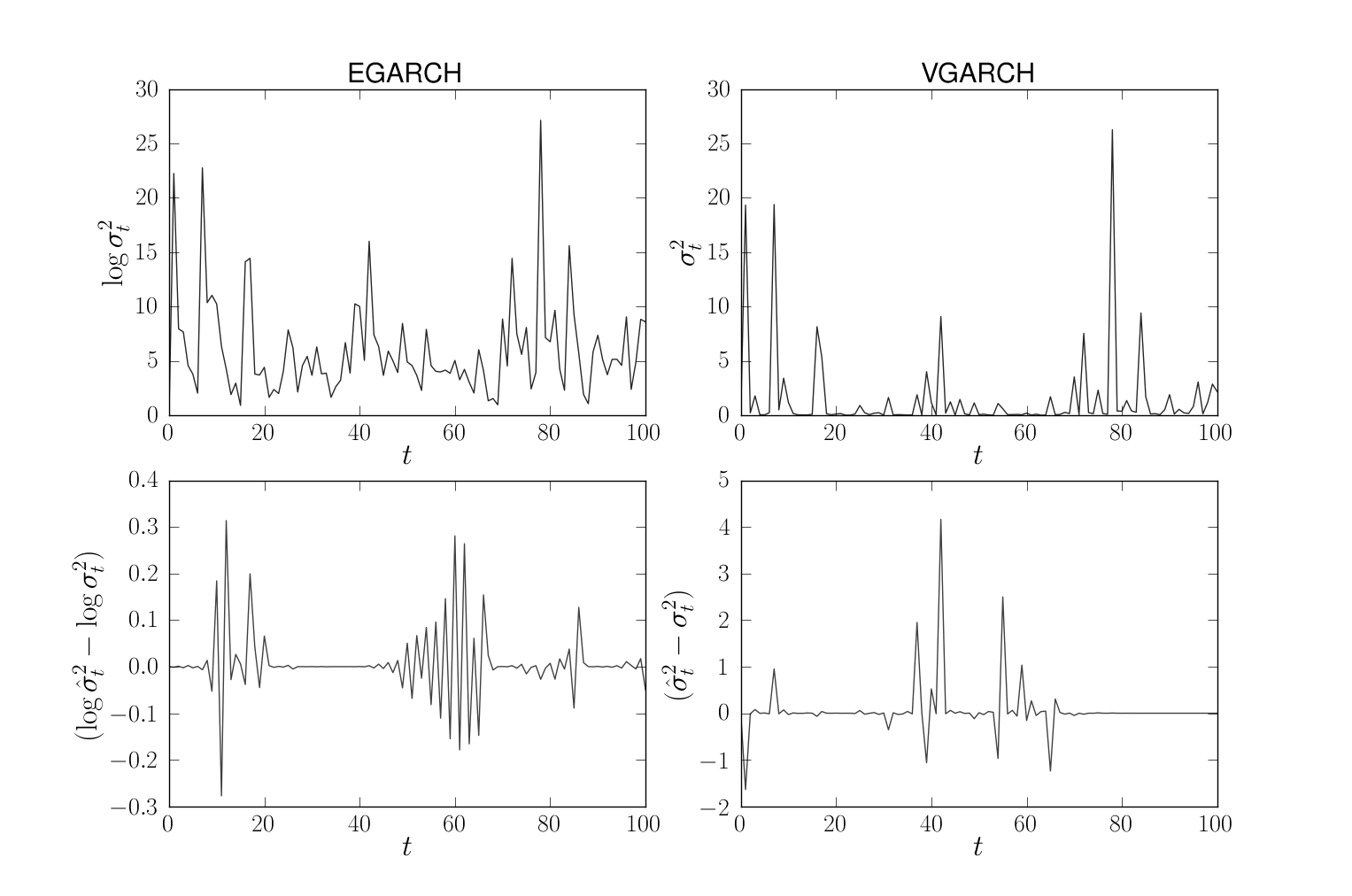}
\end{center}
\caption{EGARCH and VGARCH simulation results. The top plot shows simulated 
$\sigma^2_t$ as a function of $t$. The bottom plot shows the 
difference $(\hat{\sigma}^2_t-\sigma^2_t)$ as a function of $t$.}
\label{fig_sim}
\end{figure}

Finally, it is possible in principle that the conditions ensuring the
existence of the stationary solution and the lack of stability are mutually
exclusive and therefore the Theorems \ref{thm_stat_distr_egarch} -
\ref{thm_lln_noninv_vgarch} are vacuous. However, it is easy to check
the opposite (see also the Figures \ref{fig_lambda_y} and \ref{fig_sim}).

For EGARCH, we'll show that for any distribution of $\{\varepsilon_0\}$
satisfying Condition \ref{cond_regul} there exists $\theta$, such that
the Condition \ref{cond_stat_egarch} is
satisfied and $\lambda^y_E > 0$. Indeed, the Condition \ref{cond_stat_egarch}
is satisfied if $\beta < 1$ and $E |\log | \varepsilon_0| | < \infty$.
Condition \ref{cond_regul} ensures that the latter holds.
It therefore suffices to choose arbitrary $\alpha>0$, $\beta < 1$, $\delta = 0$ and large enough $\gamma$.

For VGARCH we'll only check a more relaxed statement. Assume in addition to
Condition \ref{cond_regul} that there exists $\mu \in \mathbb{R}$, such that 
$$
E \log |\varepsilon_0 - \mu| < E \log |\varepsilon_0|.
$$
As a particular example of a distribution with such property, it suffices
to choose a unimodal non-symmetric distribution with a maximum not at $0$.
For example, an appropriately scaled mixture of $\xi_1$ 
and $-\xi_2$, where both $\xi_1$ and $\xi_2$ had exponential distributions 
but with different parameters, is one such distribution.

Put $\beta = 0$ and select
arbitrary $\gamma > 0$, $\alpha > 0$ and $\delta = \mu$. Then
$$
\overline{\sigma}^2_0 = \alpha + \gamma (\varepsilon_{-1} - \mu)^2, \quad
\left (H^V \right )'_x(\theta, \overline{\sigma}_0 \varepsilon_0,
\overline{\sigma}^2_0) =
\frac{\gamma \varepsilon_0 (\varepsilon_0 - \mu)}
{\alpha + \gamma (\varepsilon_{-1} - \mu)^2}.
$$
Due to dominated convergence theorem, when $\alpha \to 0+$,
$$
E \log \left | \beta - \gamma \varepsilon_0 \overline{\sigma}^{-2}_0
(\varepsilon_0 - \mu) \right | =
E \log \left | \frac{\gamma \varepsilon_0 (\varepsilon_0 - \mu)}
{\alpha + \gamma (\varepsilon_{-1} - \mu)^2} \right| \to
$$
$$
E \log \left | \frac{\varepsilon_0}{\varepsilon_{-1} - \mu} \right| =
E \log | \varepsilon_0 | - E \log | \varepsilon_{-1} - \mu | > 0.
$$
For arbitrary $\gamma$ and $\alpha$ and $\beta$ small enough,
$\lambda^y_V > 0$.

To illustrate the main results, we show the figures demonstrating the
behaviour of instability coefficient $\lambda^y$ and sample paths of models
(\ref{equ_egarch}) and (\ref{equ_vgarch}). 

\subsection{Deterministic Chaos in EGARCH}\label{sec_egarch_chaos}
We'll discuss a simple case of
deterministic EGARCH model, which sheds light onto the origin of
non-invertibility in this particular model and also on the intuitive sense
of local invertibility. Assume that the Condition \ref{cond_stat_egarch}
is satisfied and
\begin{equation}\label{equ_cond_determ_egarch}
P(\varepsilon_1 = 1) = P(\varepsilon_1 = -1) = 1 / 2, \quad \beta < 1, \quad \delta = 0,
\end{equation}
then
$$
\log \overline{\sigma}_t^2 \equiv \log  \overline{\sigma}^2 =
(\alpha + \gamma) / (1 - \beta),
$$
\begin{equation}\label{equ_determ_egarch}
d_t = f(d_{t-1}), \quad f(x) = \alpha + \beta x + \gamma e^{-d_{t-1} / 2},
\quad d_t = \log \hat{\sigma}_t^2 - \log \sigma_t^2, \quad t \in \mathbb{Z}.
\end{equation}
We don't pursue the task of the formal proof, but numerical simulations show
that the dynamical system given by $f$ of (\ref{equ_determ_egarch})
is chaotic in the
sense of Devaney once $\gamma$ is large enough.
More precisely,
when $\gamma$ increases, the system
experiences period-doubling bifurcations and eventually becomes chaotic.

\begin{figure}[h]
\label{fig_determ_egarch_bif_diag}
\includegraphics[scale=0.3]{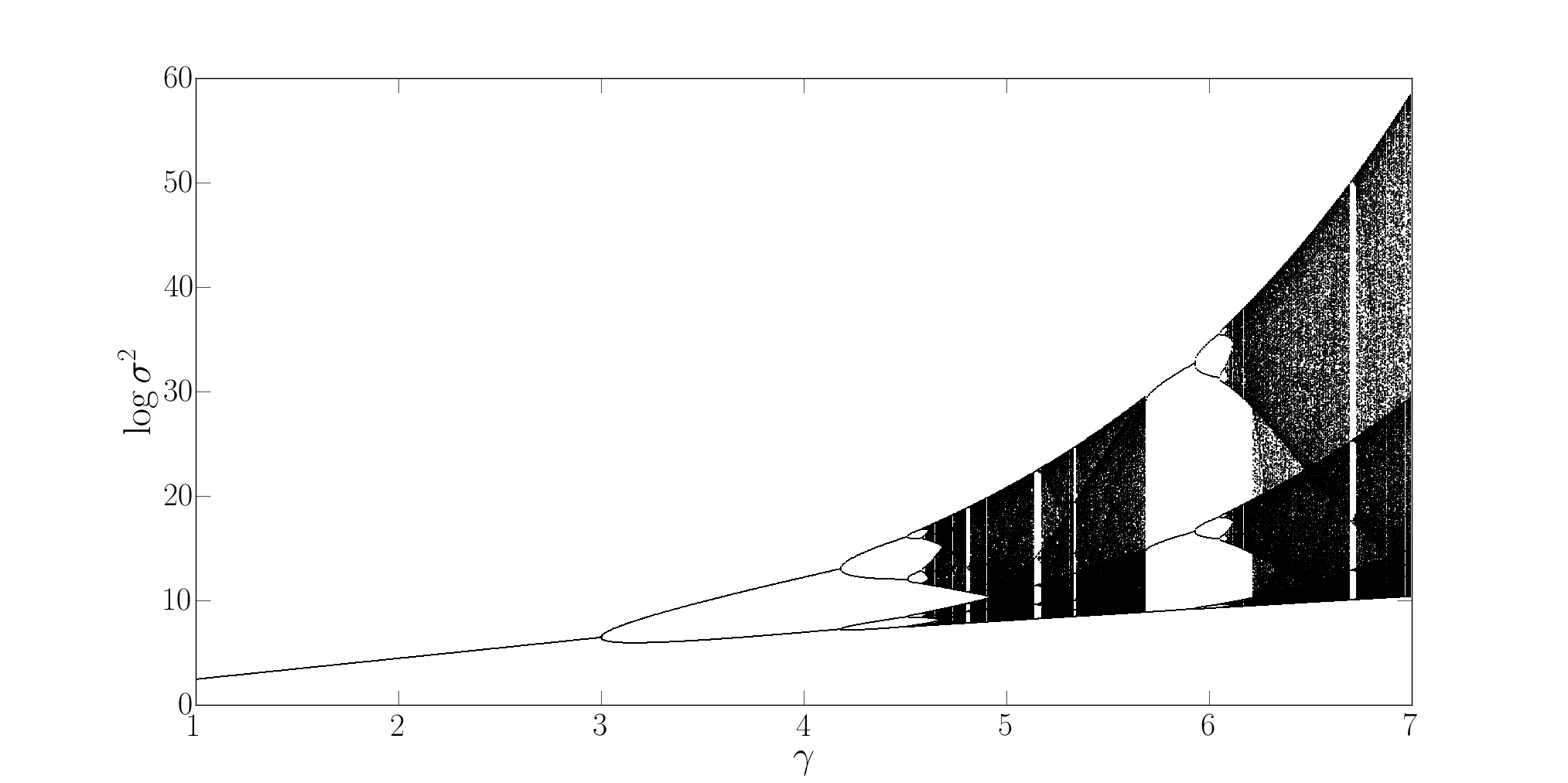}
\caption{Deterministic EGARCH bifurcation diagram for 
$\alpha = 0.2$, $\beta = 0.5$.}
\end{figure}

The bifurcation diagram of (\ref{equ_determ_egarch}) is shown at Figure
\ref{fig_determ_egarch_bif_diag}.
As discussed in \cite{Devaney1989}, the period-doubling behaviour
is typical for maps
$\mathbb{R} \to \mathbb{R}$. For example, the famous logistic map
$x \mapsto \gamma x (1 - x)$ exhibits the same type of bifurcations. It is also quite
interesting that under (\ref{equ_cond_determ_egarch}) the
local non-invertibility of EGARCH is equivalent to instability
of $\log  \overline{\sigma}^2$ which is the only
fixed point of $f$. If this point was stable in the sense
of $\lambda^y_E < 0$, there would be at least an interval of stable behaviour.
In the general "stochastic" case, discussed in the previous subsection,
the stability is equivalent to invertibility with the exception of the case 
$\lambda^y = 0$.


\section{Proofs}\label{sec_proofs}

\subsection{Markov chain foundations}
In this section, we'll present the proofs of Theorems \ref{thm_stat_distr_egarch} -
\ref{thm_lln_noninv_vgarch}. The main technical tool used is the general
result stating the existence of a stationary distribution in Markov chains
(Theorem 10.0.1 of \cite{MeynTweedie2009}), we'll re-state it for convenience.
 We also re-use some notation and definitions
from the book. 

Assume that $\{\Phi_t\}_{t \in \mathbb{Z}}$ is a Markov chain
with a polish space $\mathcal{X}$
as state space, transition probabilitites
$P(x, A), \; x \in \mathcal{X}, \; A \in \mathcal{B} (\mathcal{X})$,
$\phi$ is a measure on $\mathcal{B}(\mathcal{X})$. Also
fix $\{a_i\}$ - a distribution on $\mathbb{N}^+$. Denote
$$
\tau_A = \inf\{t \geq 0: \Phi_t \in A\}, \quad \eta_A =
\sum_{t \geq 1} I\{\Phi_t \in A\},
$$
$$
P_x(\cdot) = P(\cdot | \Phi_0 = x), \quad E_x[\cdot] = E[\cdot | \Phi_0 = x].
$$
$$
L(x, A) = P_x(\tau_A < \infty).
$$
Denote $\psi$ the maximal  irreducibility measure (see \cite[Proposition 4.2.2]{MeynTweedie2009}),
$$
\mathcal{B}^+(\mathcal{X}) = \{B \in \mathcal{B}(\mathcal{X}): \psi (B) > 0\}.
$$

\begin{theorem}[Theorem 10.0.1 of \cite{MeynTweedie2009}]
\label{thm_inv_measure_existence}
If the chain $\Phi$ is recurrent then it admits a unique (up to constant
multiples) invariant measure $\pi$. This invariant measure is finite
(rather than merely $\sigma$-finite, and therefore the chain is
positive), if there exists a petite set $C$ such that
\begin{equation}\label{equ_strict_recurrence}
\sup_{x \in C} E_x [\tau_C ] < \infty.
\end{equation}
\end{theorem}


\subsection{Proofs of Theorems \ref{thm_stat_distr_egarch} and
\ref{thm_stat_distr_vgarch}}
The proofs of these two theorems will follow exactly the same recipe. In Step 1,
we establish the following drift criteria from \cite{MeynTweedie2009} for 
a certain set $C$:
\begin{definition}[Strict drift towards C]
For some set $C \in \mathcal{B}(\mathcal{X})$, some constants $b < \infty$
and $\mu > 0$
and an extended real-valued function $V : \mathcal{X} \to [0, \infty]$
\begin{equation}\label{equ_strict_drift}
\Delta V (x) \leq -\mu I\{x \notin C\} + b I\{x \in C\}, \quad x \in \mathcal{X},
\end{equation}
where
$$
\Delta V (x) = E_x [V (\Phi_1 )] - V(x), \quad  x \in \mathcal{X}.
$$
\end{definition}
In Step 2 we'll check that the chosen $C$ is petite.

In Step 3 we check the condition (\ref{equ_strict_recurrence}) and 
$\psi$-irreducibility of $\{\Phi_t\}$, using the strict drift condition. 

Finally, we establish that the appropriate Markov chain is indeed recurrent
using the following slightly modified version of Theorem 8.4.3
from \cite{MeynTweedie2009}:
\begin{theorem}\label{thm_drift_criteria}
Assume that $\{\Phi_t\}$ is $\psi$-irreducible and there exist a petite
set $C$ and a function $V$ which is unbounded off
petite sets (in the sense that $C_V (n) = \{y: V(y) \leq n\}$ is
petite for all n), such that $\Phi$ has a strict drift towards $C$. Then
$L(x, C) = 1$ for all $x \in \mathcal{X}$ and $\Phi$ is recurrent.
\end{theorem}
The recurrence of $\{\Phi_t\}$ and Theorem \ref{thm_inv_measure_existence}
imply the existence of the stationary distribution $\pi$, and finally
the equation (\ref{equ_strict_recurrence}) ensures $\pi$ is finite.

Throughout the proofs we will denote
$o_{B}(1)$ ($O_{B}(1)$) a function of $B$ which converges to $0$
(is bounded) as $B \to \infty$ and model parameters and the
distribution of $\varepsilon_0$ are fixed.


\textbf{Proof of Theorem \ref{thm_stat_distr_egarch}.}

First of all, note that the equation (\ref{equ_egarch_sigmahat})
together with (\ref{equ_egarch}) are equivalent to
\begin{equation}\label{equ_egarch_evol_zhat}
\begin{cases}
z_t = \alpha + \beta z_{t-1} + \gamma |\varepsilon_{t-1}| + \delta \varepsilon_{t-1} ,\\
\hat{z}_t = \alpha + \beta \hat{z}_{t-1} + (\gamma |\varepsilon_{t-1}| + \delta \varepsilon_{t-1})
 \exp\{-d_{t - 1} / 2\},
\end{cases}
\quad t \geq 1,
\end{equation}
where we define
$$
z_t = \log \sigma^2_t, \quad \hat{z}_t = \log \hat{\sigma}^2_t(\theta),
\quad d_t = \hat{z}_t - z_t, \quad t \geq 0.
$$
It is more convenient to instead consider the chain $\{(z_t, d_t)\}$ defined by
\begin{equation}\label{equ_egarch_evol}
\begin{cases}
z_t = \alpha + \beta z_{t-1} + \varepsilon^*_{t-1},\\
d_t = \beta d_{t-1} + \varepsilon^*_{t-1} (\exp\{-d_{t - 1} / 2\} - 1),
\end{cases}
\quad t \geq 1,
\end{equation}
where for brevity we put
$$
\varepsilon^*_t = \gamma |\varepsilon_{t}| + \delta \varepsilon_{t}, \quad t \geq 0,
$$
denote for $x \in \mathbb{R}$  
$$
e_x = \exp\{-x/2\} - 1, \quad x^+ = \max(x, 0)
$$
and recursively define 
$$
h_1(z, x) = \alpha + \beta z + x, \quad h_t(z, x_1, \ldots, x_t) =
 h_1(h_{t-1}(z, x_1, \ldots, x_{t-1}), x_t), \quad t \geq 2,
$$
$$
\hat{h}_1(d, x) = \beta d + x  e_d, \quad \hat{h}_t(d, x_1, \ldots, x_t) =
\hat{h}_1(\hat{h}_{t-1}(d, x_1, \ldots, x_{t-1}), x_t), \quad t \geq 2,
$$
so that
$$
z_t = h_t(z_0, \varepsilon^*_0, \ldots, \varepsilon^*_{t-1}), \quad 
d_t = \hat{h}_t(d_0, \varepsilon^*_0, \ldots, \varepsilon^*_{t-1}).
$$

It is easy to see that for any $d \neq 0$ and $t \geq 0$ under Condition
\ref{cond_regul} $P_{(z, d)} (d_t = 0) = 0$. Indeed, due to Condition \ref{cond_regul}
$P(\varepsilon^*_{0} = x) = 0$ for any $x \neq 0$ and 
$$
P_{(z, d)} (d_1 = 0) = P_{(z, d)} \Bigl ( \beta d_0 + \varepsilon^*_0 e_d = 0 \Bigr )
= P_{(z, d)} \Bigl (  \varepsilon^*_{0} = - \beta d_0 e_d^{-1} \Bigr ) = 0.
$$
For an arbitrary $t \geq 0$ we get the result by induction and Fubini theorem. 
From now on we will consider the modification of the chain  $\{(z_t, d_t)\}$ on the space 
$$
\mathcal{X}_E = \mathbb{R} \times (\mathbb{R} \backslash \{0\})
$$
(which exists as we've just shown).

\textbf{Step 1.} 
Put for some $k > 0$
$$
V_E(z, d) = R(d) + k |z|, \quad R(d) = \begin{cases} - d - 1, & \mbox{if } d < -1 \\
- \log |d|, & \mbox{if } |d| \leq 1 \\ \log d, & \mbox{if } d > 1 \end{cases},
$$
$$
C_E(A, B) = \left \{(z,d) \in \mathcal{X}_E:
 0 \leq |z| \leq A, |d| \geq B^{-1}, |d| \leq B \right \}.
$$
We'll first show that for some $\mu > 0$
and all $(z, d) \in \mathcal{X}_E \backslash C_E(A, B)$
\begin{equation}\label{equ_VE_estimate}
E_{(z, d)} [V_E(z_1, d_1)] \leq V_E(z, d) - \mu
\end{equation}
if $A$ and $B$ are large enough. Indeed,
\begin{equation}\label{equ_RE_estimate_z}
E_{z} |z_1| \leq |\alpha| + \beta |z| + E |\varepsilon^*_0| \leq
O_B(1) + \beta |z|.
\end{equation}

Note that the conditions $\lambda^y_E > 0$ and $\gamma \geq |\delta|$ ensure that $\gamma > 0$.
For $d \geq B$ and $B \geq \beta^{-1}$, by the definition
of $R$ and Lemma \ref{lem_log_estimate} (\ref{equ_log_estimate_b})
$$
E_{d} R(d_1) = E \log(\beta d + \varepsilon^*_0 e_d)
I \left \{d_1 \geq 1 \right\} +
$$
$$
E \bigl [ - \log \left|\beta d + \varepsilon^*_0 e_d \right|
\bigr ] I \left \{ |d_1| < 1 \right \} +
$$
$$
E \left (- \beta d - \varepsilon^*_0 e_d - 1
\right ) I \left \{d_1 \leq - 1 \right \} \leq
$$
\begin{equation}\label{equ_RE_estimate_d_geq_B}
P (d_1 \geq 1) \log (\beta d) + o_B(1) +
O_B(1) E \varepsilon^*_0 P(d_1 < - 1)
\leq \log d + \log \beta + o_B(1).
\end{equation}
For $d \leq - B$, once again using definition of $R$ and Lemma
\ref{lem_log_estimate} (\ref{equ_log_estimate_c}), we get
$$
E_d R(d_1) = E \bigl (- \beta d - \varepsilon^*_0 e_d - 1 \bigr)
I \left \{d_1 \leq - 1 \right\} +
$$
$$
E \Bigl [ - \log \bigl |\beta d + \varepsilon^*_0 e_d \bigr | \Bigr]
I \left \{|d_1| < 1\right \} +
$$
$$
E \log \bigl (\beta d + \varepsilon^*_0 e_d \bigr )
I \left \{d_1 \geq 1\right\} \leq
$$
$$
- \beta d  P \left \{d_1 \leq - 1 \right\} +
\log (- \beta d) P (|d_1| < 1) + \log (\exp\{-d/2\}) +
 E \log^+ (\varepsilon^*_1) +  O_B(1) \leq
$$
\begin{equation}\label{equ_RE_estimate_d_leq_minus_B}
- d /2 + O_B(1).
\end{equation}

For $d \in (0, B^{-1})$, using Taylor's expansion we find that for any
$\rho > 0$
$$
E_d R(d_1) \leq E \bigl (- \beta d - \varepsilon_0 e_d - 1 \bigr)
 I \left \{d_1 \leq - 1\right \} +
$$
$$
E \Bigl | \log \bigl |\beta d + \varepsilon^*_0 e_d \bigr | \Bigr| \leq
$$
\begin{equation}\label{equ_RE_estimate_d_leq_B_inv}
\leq (-e_d) E \varepsilon^*_0 - \log d + \sup_{x \in (1-\rho, 1+\rho)}
E \left |\log \left |\beta - x 2^{-1} \varepsilon^*_0 \right | \right|
+o_B(1).
\end{equation}
Using condition $\lambda^y_E > 0$
and the dominated convergence theorem, we infer that  for some $\mu > 0$
and $B$ large enough the right-hand
side of (\ref{equ_RE_estimate_d_leq_B_inv}) is no more than
$$
-\log d - \mu, \quad \mu > 0.
$$
The case  $d \in (-B^{-1}, 0)$ is similar.

It is now enough to notice that equations (\ref{equ_RE_estimate_z}) -
(\ref{equ_RE_estimate_d_leq_B_inv}) indeed imply (\ref{equ_VE_estimate}).

\textbf{Step 2.}
We'll show now that for some $\mu > 0$, $D$ the sets $C_E(A, B)$
are $\nu$-petite for the chain $\{(z_t, d_t)\}$, $a = (0, 1/2, 1/2, 0, \ldots)$ and
$$
d \nu = \mu I\{(z, d) \in Z_E\} d \lambda_E, 
\quad Z_E = [D, D + \mu] \times [D, D + \mu],
$$
where $\lambda_E$ is the Lebesgue measure on $\mathcal{X}_E$.
It will also ensure that $V_E$ is
unbounded off petite sets, since it easy to see
that sets of form $\{V_E(z, d) \leq x\}$ are subsets of
$C_E(A, B)$ for large enough $A$, $B$.

First, consider the case $d > 0$. We'll use the Lemma
\ref{lem_jacobian} with
$$
k = 2, \quad \tau = (z, d), \mathcal{T} = C_E^+(A, B) = C_E(A, B) \cap \{d > 0\}, \quad
Z = Z_E, 
$$
$$
f_{(z, d)}(x_0, x_1) = (h_2(z, x_0, x_1), \hat{h}_2(d, x_0, x_1)).
$$

The only non-trivial part is checking that once $D$ is large enough
the equation
\begin{equation}\label{equ_petiteness_equation_egarch}
h_2(z, x_0, x_1) = z^*, \quad \hat{h}_2(d, x_0, x_1) = d^*
\end{equation}
on $x_0$ and $x_1$ has a solution for any
$(z^*, d^*) \in Z_E$ and $(z, d) \in C_E^+(A, B)$.  
The $z$ part of (\ref{equ_petiteness_equation_egarch}) is equivalent to
$$
z^* = \alpha + \beta (\alpha + \beta z + x_0) + x_1,
$$
which is satisfied if and only if
$$
x_0 = \beta^{-1} (S - x_1).
$$
where we put
$$
S = z^* - \alpha - \beta (\alpha + \beta z).
$$
Then for $x_0 = S / (2 \beta)$, $x_1 = S / 2$ the $z$ part of
(\ref{equ_petiteness_equation_egarch}) is satisfied and for large enough $D$ we have
$$
\hat{h}_2 (d, S / (2 \beta), S / 2) = \beta \left (\beta d + x_0 e_d \right)
+ x_1 \left ( e^{-(\beta d + x_0 e_d) / 2} - 1 \right) \geq
$$
\begin{equation}\label{equ_petiteness_equation_egarch_1}
\beta^2 d + S e_d / 2 + (S/2) \left ( e^{-(\beta d + e_d S / (2 \beta)) / 2} - 1 \right)
\geq (S/2) \left (e^{- e_{-1/(2 B)} S / 4 \beta - \beta B / 2} - 2 \right ) 
\geq D,
\end{equation}
where we used the definition of $C_E(A,B)$. On the other hand, if
$x_0 = \beta$, $x_1 = S - \beta^2$ and $D$ is large enough
$$
\hat{h}_2 (d, 0, S) = \beta \left (\beta d + x_0 e_d \right)
+ x_1 \left ( e^{-(\beta d + x_0 e_d) / 2} - 1 \right) =
$$
\begin{equation}\label{equ_petiteness_equation_egarch_2}
\beta^2 d + (S - \beta^2) \left ( e^{-\beta d / 4} - 1 \right) \leq 
\beta^2 d + (S - \beta^2) \left ( e^{-\beta B^{-1} / 4} - 1 \right) \leq 0. 
\end{equation}
From (\ref{equ_petiteness_equation_egarch_1}) 
and (\ref{equ_petiteness_equation_egarch_2}) and continuity of $\hat{h}_2$ 
we conclude that for some $x_1 = x_1^* \in [S / 2, S]$ and $x_0 = \beta^{-1} (S - x^*_1)$
(\ref{equ_petiteness_equation_egarch}) is satisfied.
Lemma \ref{lem_jacobian} implies that $C^+_E(A, B)$ is petite. The fact that
$C_E^-(A, B) = C_E(A, B) \cap \{d < 0\}$ is petite is proved by using
Lemma \ref{lem_jacobian_more} with
$$
k = 2, \quad \tau_1 = (z, d), \quad \tau_2 = , \quad \mathcal{T}_1 = C^-_E(A, B), \quad
Z = Z_E, \quad \mathcal{T}_2 = \mathbb{R}, \quad \mathcal{T} = C^-_E(A, B) \times \mathbb{R}^+,
$$
$$
\boldsymbol{\xi} = (\varepsilon^*_1, \varepsilon^*_2), \quad
\chi = \varepsilon^*_0, \quad f_{((z, d), x_0)}(x_1, x_2) = 
(h_3(z, x_0, x_1, x_2), \hat{h}_3(d, x_0, x_1, x_2)).
$$
The technical details are omitted as they are similar to the above case.
The set $C_E(A,B)$ is petite as a union of 
$C^+_E(A, B)$ and $C^-_E(A, B)$.

In order to ensure the strict drift condition (\ref{equ_strict_drift}),
it only remains to show that
$$
\sup_{z \in \mathbb{R}, d \in \mathbb{R} \backslash \{0\}}
(E_{(z,d)} [ V_E(z_1, d_1)] - V_E(z, d)) < \infty.
$$
For $z$ part of $V$,
$$
\sup_{z \in \mathbb{R}} \left (E_z |z_1| - |z| \right ) =
\sup_{z \in \mathbb{R}} \left (\alpha + \beta |z| +
E |\varepsilon^*_0| - |z| \right ) < \infty.
$$
For $d$ part, it is achieved in a similar fashion to the proofs of
(\ref{equ_RE_estimate_z}) -
(\ref{equ_RE_estimate_d_leq_minus_B}). Thus, the correctness of
(\ref{equ_strict_drift}) is established. Theorem \ref{thm_drift_criteria}
implies that the chain $\{(z_t, d_t)\}$ on the space $\mathcal{X}_E$
is indeed recurrent.

\textbf{Step 3.} Let's check now that (\ref{equ_strict_recurrence}) takes
place with $\Phi_k = (z_k, d_k)$, $C = C_E(A, B)$ and $A$ and $B$ large enough. Put
$$
Q_k = \{(z_k, d_k ) \in C_E(A, B), \; (z_i, d_i ) \notin C_E(A,B), \;
i = 1, \ldots, k -1\}, \quad x = (z, d)
$$
and note that due to (\ref{equ_VE_estimate}) and law of iterated expectations
$$
\infty > E_x [V_E(\Phi_1)] \geq E_x [V_E(\Phi_1) I\{\Phi_1 \notin C\}] \geq
E_x [(V_E(\Phi_2) + \mu) I\{\Phi_1 \notin C\}] \geq
$$
$$
\mu P_x(\Phi_1 \notin C) + E_x [V_E(\Phi_2)
I\{\Phi_1 \notin C, \Phi_2 \notin C\}] \geq
$$
$$
\mu P_x(\Phi_1 \notin C) + \mu P_x(\Phi_2 \notin C,
\Phi_1 \notin C) + E_x [V_E(\Phi_3)
I\{\Phi_1 \notin C, \Phi_2 \notin C, \Phi_3 \notin C\}] \geq \ldots
$$
$$
\mu \sum_{k \geq 1} P_x(\Phi_1 \notin C, \ldots \Phi_k \notin C) \geq
\mu \sum_{k \geq 1} k P_x(\Phi_1 \notin C, \ldots \Phi_k \notin C,
 \Phi_{k+1} \in C) \geq
$$
$$
\mu (E_x [\tau_C] - 1).
$$
We'll now prove that for some $D$, $\mu > 0$
 the chain $\Phi$ is $\phi$-irreducible with 
$$
d \phi = \mu I\{Z_E\} d \lambda_E.
$$
Indeed, for any $x \in C_E(A, B)$ and $G \subseteq Z_E$
$$
P_x(\Phi_2 \in G) \geq \mu \lambda_E(G), 
$$
hence for any $x \in C_E(A, B)$ and 
$G \subseteq Z_E$
$$
L(x, G) \geq \phi(G).
$$
For any $x \in \mathcal{X}_E \backslash C_E(A, B)$ due to (\ref{equ_strict_recurrence})
$$
L(x, C_E(A, B)) = 1
$$
and hence there exists $k$, $P_x(\Phi_k \in C_E(A, B)) > 0$. Therefore
for any $G \subseteq Z_E$, $\lambda_E(G) > 0$
$$
L(x, G) \geq P_x(\Phi_{k+2} \in G) \geq P_x(\Phi_k \in C_E(A, B), 
\Phi_{k + 2} \in G) \geq P_x(\Phi_k \in C_E(A, B)) \mu \lambda_E(G) > 0.
$$

Now all the conditions of Theorem \ref{thm_inv_measure_existence}
are verified. Hence, we have established that the system
(\ref{equ_egarch_evol}) has an invariant probability measure. Define a new  
(extended) probability space $(\Omega^*, \mathcal{F}^*, P^*)$ as 
a cartesian product
$$
\Omega^* = \Omega \times \mathbb{R}, \quad \mathcal{F}^* = 
\mathcal{F} \times \mathcal{B}(R), \quad P^* = P \times U,
$$
where $U$ is a uniform measure concentrated on $[0, 1]$. There exists
a measurable function $r: \mathbb{R}^2 \to \mathbb{R}$, such that the distribution of 
$$
(\overline{\sigma}^2_0, r(\overline{\sigma}^2_0, u))
$$
is the invariant probability measure of (\ref{equ_egarch_evol}). Put
$$
\overline{d}_0 = (\overline{\sigma}^2_0, r(\overline{\sigma}^2_0, u)), \quad 
\overline{d}_t = \hat{h}(\overline{d}_{t-1}, \varepsilon_{t-1}^*), \quad t \geq 1.
$$
Then $\{(\overline{z}_t, \overline{d}_t)\}$ is a
strictly stationary solution of (\ref{equ_egarch_evol}). It only remains to note that
due to the definition of $\mathcal{X}_E$
$$
\overline{\sigma}^2_t = \overline{z}_t, \quad
\hat{\overline{\sigma}}^2_t = \overline{z}_t + \overline{d}_t
$$
is the stationary solution of (\ref{equ_sigma_sigmahat_system}) with
$H = H_E$ and
$$
P(\overline{\sigma}_0^2 = \hat{\overline{\sigma}}_0^2) = 0. \Box
$$


\textbf{Proof of Theorem \ref{thm_stat_distr_vgarch}.}

Let's rewrite the VGARCH evolution equations (\ref{equ_vgarch})
and (\ref{equ_vgarch_sigmahat})
in the explicit form:
\begin{equation}\label{equ_vgarch_evol}
\begin{cases}
z_t = \alpha + \beta z_{t-1} + \gamma (\varepsilon_{t-1} - \delta)^2,\\
\hat{z}_t = \alpha + \beta \hat{z}_{t-1} + \gamma
\left (\varepsilon_{t-1} z_{t-1}^{1/2} \hat{z}_{t-1}^{-1/2}  - \delta \right)^2,
\end{cases}
\quad t \geq 1,
\end{equation}
where we define
$$
z_t = \sigma^2_t, \quad \hat{z}_t = \hat{\sigma}^2_t, \quad t \geq 0.
$$
We will consider a Markov chain $\{(z_t, \hat{z}_t)\}$. Define
$$
\mathcal{X}_V = \{(z, \hat{z}) \in \mathbb{R}^2: z \geq \alpha,
\hat{z} \geq \alpha, z \neq \hat{z}\}.
$$
If $(z, \hat{z}) \in \mathcal{X}_V$, then
$$
P_{(z, \hat{z})}((z_1, \hat{z}_1) \notin \mathcal{X}_V) = P_{(z, \hat{z})}\left (\alpha + \beta z + \gamma (\varepsilon_0 - \delta)^2 =  \alpha + \beta \hat{z} + \gamma
\left (\varepsilon_0 z^{1/2} \hat{z}^{-1/2}  - \delta \right)^2 \right )  =
$$
\begin{equation}\label{equ_why_chain_stays_in_xv}
P_{(z, \hat{z})} \left (\beta (z - \hat{z}) + \gamma \varepsilon_0 (\varepsilon_0 (1 - z \hat{z}^{-1}) + 2 \delta (z^{1/2} \hat{z}^{-1/2} - 1) ) = 0 \right ) = 0,
\end{equation}
because the second degree polynomial with respect to $\varepsilon_0$ 
inside the probability in (\ref{equ_why_chain_stays_in_xv}) has at most 
$2$ roots and the distribution of $\varepsilon_0$ is continuous due
 to Condition \ref{cond_regul}. Similarly, for any 
 $(z, \hat{z}) \in \mathcal{X}_V$ and $t > 0$ due to induction and Fubini theorem
$$
P_{(z, \hat{z})}((z_t, \hat{z}_t) \notin \mathcal{X}_V) = 0.
$$
Hence we the chain $\{(z_t, \hat{z}_t)\}$ has a modification defined on $ \mathcal{X}_V$, which we will consider from now on.

We also recursively define 
$$
h_1(z, x) = \alpha + \beta z + \gamma (x - \delta)^2,
 \quad h_t(z, x_1, \ldots, x_t) =
 h_1(h_{t-1}(z, x_1, \ldots, x_{t-1}), x_t), \quad t \geq 2,
$$
$$
\hat{h}_1(z, \hat{z}, x) = \alpha + \beta \hat{z} + \gamma 
(x z^{1/2} \hat{z}^{-1/2} - \delta)^2, \quad \hat{h}_t(z, \hat{z}, x_1, \ldots, x_t) =
\hat{h}_1(\hat{h}_{t-1}(z, \hat{z}, x_1, \ldots, x_{t-1}), x_t), \quad t \geq 2,
$$
so that
$$
z_t = h_t(z_0, \varepsilon_0, \ldots, \varepsilon_{t-1}), \quad 
\hat{z}_t = \hat{h}_t(z_0, \hat{z}_0, \varepsilon_0, \ldots, \varepsilon_{t-1}).
$$

\textbf{Step 1.} Fix positive real $k, m, A, B$ and $n \in \mathbb{N}$.
Their values will be chosen later, for now we'll just say that we select 
them in the following order: $k$ is small enough, then $n$ large enough, 
$B$ then large enough, $m$ large enough, and finally $A$ - large enough. Define
$$
R_n(z, \hat{z}) = E_{(z, \hat{z})} [- \log (|z_n - \hat{z}_n|)],
\quad V_V(z, \hat{z}) = z + k \hat{z} + m R_n(z, \hat{z}), \quad z \neq \hat{z},
$$
$$
C_V(A, B) = \{(z, \hat{z}) \in \mathbb{R}^2: z \in [\alpha, A],
\hat{z} \in [\alpha, A], |z - \hat{z}| \geq B^{-1} \}.
$$
In order to check the strict drift condition, we'll demostrate that for some
$\mu > 0$ and all $(z, \hat{z}) \in \mathcal{X}_V \backslash C_V(A, B)$
\begin{equation}\label{equ_VV_estimate}
E_{(z, \hat{z})}  [V_V(z_1, \hat{z}_1)] \leq V_V(z, \hat{z}) - \mu
\end{equation}
if $k, m, A, B$ and $n$ are chosen appropriately. First note that
\begin{equation}\label{equ_RV_estimate_simple_z_zhat}
\begin{split}
E_{(z, \hat{z})} z_1 = \alpha + \beta z + \gamma (\delta^2 + 1),\\
E_{(z, \hat{z})} \hat{z}_1 = \alpha + \beta \hat{z} + \gamma
\left (\delta^2 + z \hat{z}^{-1} \right) \leq \alpha + \beta \hat{z}
+ \gamma \delta^2 + \alpha^{-1} \gamma z.
\end{split}
\end{equation}
Second, we'll show that for any $(z, \hat{z}) \in \mathcal{X}_V$
\begin{equation}\label{equ_RV_estimate_any_z_log_diff}
E_{(z, \hat{z})}[\log(|z - \hat{z}|) - \log(|\hat{z}_1 - z_1|)] \leq M,
\end{equation}
where $M$ depends only on model parameters. 
For any $(z, \hat{z}) \in \mathcal{X}_V$
$$
\hat{z}_1 - z_1 = \beta (\hat{z_0} - z_0) + \gamma \left ( \varepsilon_0^2
\left ( \frac{z_0}{\hat{z_0}} - 1\right ) - 2 \varepsilon_0 \delta
\left ( \frac{z_0^{1/2}}{\hat{z_0}^{1/2}} - 1\right ) \right ),
$$
\begin{equation}\label{equ_RV_estimate_any_z_log_diff_interm}
[ - \log(|\hat{z}_1 - z_1|)] = - \log(\beta) - \log(|\hat{z} - z|) +
\left [ -\log (|1 + q \varepsilon_0 - r \varepsilon_0^2 |) \right ],
\end{equation}
where
$$
q = \frac{2 \gamma \delta }
{\beta \hat{z}^{1/2} \left (z^{1/2} + \hat{z}^{-1/2} \right )}, 
\qquad r = \frac{\gamma }{\beta \hat{z} }.
$$
Condition $\lambda^y_V > 0$ implies that $\gamma > 0$, thus $r > 0$.
Note also that both $q$ and $r$ are uniformly bounded.
Due to Lemma \ref{lem_log_estimate} (\ref{equ_log_estimate_f}) and
(\ref{equ_RV_estimate_any_z_log_diff_interm}),
(\ref{equ_RV_estimate_any_z_log_diff}) is fulfilled.
Applying it together with
(\ref{equ_RV_estimate_simple_z_zhat})
iteratively, we obtain 
\begin{equation}\label{equ_RV_estimate_any_z_log_diff_n}
|R_n(z, \hat{z}) +  \log(|\hat{z} - z|) | \leq 
n M^*,
\end{equation}
where (not importantly), 
\begin{equation}\label{equ_M_star}
M^* = M + \alpha + \gamma (\delta^2 + \gamma \alpha^{-1} (\alpha +
\gamma (1 + \delta^2))) < \infty.
\end{equation}
For $\hat{z} > A$, $\hat{z} \geq z$,
using (\ref{equ_RV_estimate_any_z_log_diff}) and
(\ref{equ_RV_estimate_simple_z_zhat}) we infer
$$
E_{(z, \hat{z})}  [V_V(z_1, \hat{z}_1)] = E_{(z, \hat{z})} z_1 +
k E_{(z, \hat{z})} \hat{z}_1 - m E_{(z, \hat{z})} \log(|\hat{z}_n - z_n|) \leq
$$
$$
\left (\alpha + \beta z + \gamma (1 + \delta^2) \right) +
k \left (\alpha + \beta \hat{z} + \gamma (\alpha^{-1} z +
\delta^2) \right) - m \log (|\hat{z} - z|) + m n M^* \leq
$$
$$
V_V(z, \hat{z}) + (\beta - 1 + k \alpha^{-1} \gamma) z + k (\beta - 1) \hat{z} + 
\alpha + \gamma (\delta^2 + 1) + k (\alpha + \gamma \delta^2) + m n M^* \leq V_V(z, \hat{z}) - \mu
$$
for any $\mu > 0$ and $k < (1 - \beta) \alpha / \gamma$, 
\begin{equation}\label{equ_RV_estimate_zhat_geq_A}
A \geq \frac{\mu + \alpha + \gamma (\delta^2 + 1) + k (\alpha + \gamma \delta^2) 
+  m n M^*}{k (1 - \beta)}.
\end{equation}
If $z > A$, $z \geq \hat{z}$, we similarly have
$$
E_{(z, \hat{z})}  [V_V(z_1, \hat{z}_1)] \leq 
V_V(z, \hat{z}) + (\beta - 1 + k \alpha^{-1} \gamma) z + k (\beta - 1) \hat{z} + 
\alpha + \gamma (\delta^2 + 1) + k (\alpha + \gamma \delta^2) + 
$$
$$
m n M^* \leq V_V(z, \hat{z}) - \mu
$$
once 
\begin{equation}\label{equ_RV_estimate_z_geq_A}
A \geq \frac{\mu + \alpha + \gamma (\delta^2 + 1) + k (\alpha + \gamma \delta^2) 
+  m n M^*}{k (1 - \beta - k \alpha^{-1} \gamma)}.
\end{equation}
The case $z \leq A, \hat{z} \leq A$ is less straightforward.
First, fix any $n \in \mathbb{N}$ and assume $A \geq e^{\sqrt{n}}$
(any function which grows quicker than linearly but slower than
exponentially will work). The case of $(z, \hat{z})$ such that
$$
|z - \hat{z}| \leq B^{-1}, \quad z > e^{\sqrt{n}},
$$
is considered similarly, we have
$$
E_{(z, \hat{z})}  [V_V(z_1, \hat{z}_1)] \leq  V_V(z, \hat{z}) - \mu,
$$
if 
\begin{equation}\label{equ_RV_estimate_z_geq_exp_sqrt_n}
e^{\sqrt{n}} \geq \frac{\mu + \alpha + \gamma (\delta^2 + 1) + k (\alpha + \gamma \delta^2) 
+  m n M^*}{k (1 - \beta - k \alpha^{-1} \gamma)}.
\end{equation}
Consider now the case $z \leq e^{\sqrt{n}} $ (and by definition of $\mathcal{X}_V$
$|\hat{z} - z| < B^{-1}$). We will
choose $B$ large enough, depending on $n$. 
Note that
$$
E_{(z, \hat{z})} R_n(z_1, \hat{z}_1) - R_n(z, \hat{z}) =
R_{n + 1}(z, \hat{z}) - R_n(z, \hat{z}) =
$$
\begin{equation}\label{equ_RV_estimate_z_leq_exp_sqrt_n}
- E_{(z, \hat{z})} \log \left |\beta - \gamma \varepsilon_n \hat{z}^{-1}_n 
\left( \varepsilon_n - 2 \delta (1 + z_n^{1/2} \hat{z}_n^{-1/2} )^{-1} \right)
\right |.
\end{equation}
Let's estimate the right-hand side of (\ref{equ_RV_estimate_z_leq_exp_sqrt_n})
from above. Note that
$$
0 \leq \hat{z}^{-1}_n \leq \alpha^{-1}, \quad 0 \leq 
\hat{z}^{-1}_n (1 + z_n^{1/2} \hat{z}_n^{-1/2} )^{-1} \leq \alpha^{-1}
$$
thus due to Lemma \ref{lem_log_estimate} (\ref{equ_log_estimate_f}),
the family of conditional distributions
$$
\left \{ P_{(z, \hat{z})} \left ( \log \left |\beta - \gamma \varepsilon_n \hat{z}^{-1}_n 
\left( \varepsilon_n - 2 \delta (1 + z_n^{1/2} \hat{z}_n^{-1/2} )^{-1} \right)
\right | \in \cdot \right ) \right \}_{(z, \hat{z}) \in C_V(A, B), n \in \mathbb{R}}
$$
is uniformly integrable. For a fixed $n$, due to the fact that the function
$\hat{h}$ is continuous 
$$
\hat{z}_n \stackrel{P} \to z_n
$$
as $B \to \infty$ uniformly in 
$(z, \hat{z}) \in (\mathcal{X}_V \backslash C_V(A, B)) \cap \{z \leq e^{\sqrt{n}}\}$.
Therefore, due to dominated 
convergence theorem for $n \to \infty$ and $B \to \infty$ appropriately quickly,
$$
E \log \left |\beta - \gamma \varepsilon_n \hat{z}^{-1}_n 
\left( \varepsilon_n - 2 \delta (1 + z_n^{1/2} \hat{z}_n^{-1/2} )^{-1} \right)
\right | \to - \lambda^y_V.
$$
Finally, for large enough $n$ and appropriately chosen $B$, for any 
$(z, \hat{z}) \in (\mathcal{X}_V \backslash C_V(A, B)) \cap \{z \leq e^{\sqrt{n}}\}$
the right-hand side of (\ref{equ_RV_estimate_z_leq_exp_sqrt_n}) is no greater than
$$
- \lambda^y_V / 2 < 0
$$
and for all
$(z, \hat{z}) \in (\mathcal{X}_V \backslash C_V(A, B)) \cap \{z \leq e^{\sqrt{n}}\}$
$$
E_{(z, \hat{z})}  [V_V(z_1, \hat{z}_1)] \leq 
V_V(z, \hat{z}) + (\beta - 1 + k \alpha^{-1} \gamma) z + k (\beta - 1) \hat{z} + 
$$
$$
\alpha + \gamma (\delta^2 + 1) + k (\alpha + \gamma \delta^2) - m \lambda^y_V / 2
\leq V_V(z, \hat{z}) - \mu,
$$

if only
\begin{equation}\label{equ_RV_estimate_z_close_to_zhat_final}
m \geq \frac{\mu + \alpha + \gamma (\delta^2 + 1) + k (\alpha + \gamma \delta^2)}{\lambda^y_V / 2},
\end{equation}
To sum up, we have shown that for
(\ref{equ_VV_estimate}) to be true the constants
defining the function $V_V(z, \hat{z})$ may be chosen in the following order:

1) $k = (1- \beta) \alpha / (2 \gamma)$.

2) $m$ is large enough so that (\ref{equ_RV_estimate_z_close_to_zhat_final})
holds.

3) $n$ is large enough and $B$ - large enough depending on $n$, 
so that (\ref{equ_RV_estimate_z_geq_exp_sqrt_n}) holds and
the rhs of (\ref{equ_RV_estimate_z_leq_exp_sqrt_n}) is no 
greater than $- \lambda^y_V / 2$ for $(z, \hat{z}) \in 
(\mathcal{X}_V \backslash C_V(A, B)) \cap \{z \leq e^{\sqrt{n}}\}$.

4) $A$ - large enough for (\ref{equ_RV_estimate_z_geq_A}) and
(\ref{equ_RV_estimate_zhat_geq_A}) to hold.

Thus, (\ref{equ_VV_estimate}) is established.

\textbf{Step 2.} Similarly to the proof of Step 2 of the Theorem
\ref{thm_stat_distr_egarch}, we will show now that the sets $C_V(A, B)$
are $\nu$-petite for the chain $\{(z_t, \hat{z}_t)\}$,
$a = (0, 1, 0, \ldots)$ and
$$
d \nu = \mu I\{(z, \hat{z}) \in Z_V\} d \lambda_V,
\quad Z_V = [D, D + \mu] \times [D, D + \mu]
$$
with some $\mu > 0$, $D$, where $\lambda_V$ is Lebesgue measure on $\mathcal{X}_V$.
It will also ensure that $V_V$ is
unbounded off petite sets, since it easy to see
that sets of form $\{V_V(z, \hat{z}) \leq x\}$ are subsets of
$C_V(A, B)$ for large enough $A$, $B$.

For that, we will use Lemma \ref{lem_jacobian} with
$$
k = 2, \quad \mathcal{T} = C_V(A, B), \quad Z = Z_V, \quad \tau = (z, \hat{z}),
 \quad f_{(z, \hat{z})}(x_1, x_2) = (h_2(z, x_1, x_2), \hat{h}_2(z, \hat{z}, x_1, x_2) ).
$$
To ensure the Lemma's conditions, we will show that once $D$ is
large enough the equation
\begin{equation}\label{equ_C_V_equation}
h_2(z, x_0, x_1) = z^*, \quad \hat{h}_2(z, \hat{z}, x_0, x_1) = \hat{z}^*
\end{equation}
on $x_0$ and $x_1$ has a solution for
$(z^*, \hat{z}^*) \in Z_V$  and any $(z, \hat{z}) \in C_V(A, B)$, 
such that $x_0 \in [\delta, (D \beta^{-1} \gamma^{-1} )^{1/2}]$.
Indeed, there exists a continuous function $x_1(x_0, z^*)$,
defined for large $z^*$ and any $x_0 \in [\delta, (D \beta^{-1} \gamma^{-1} )^{1/2}]$, such that
\begin{equation}\label{equ_C_V_solution}
h_2(z, x_0, x_1(x_0, z^*)) \equiv z^*.
\end{equation}
Note that on $C_V(A,B)$ $(\hat{z} / z)$ is
uniformly separated from both $0$ and $\infty$. Assume that $\hat{z} > z$;
the case $\hat{z} < z$ is similar. When $x_0 =
\delta$, $\hat{h}_1(z, \hat{z}, x_0) > h_1(z, x_0)$ and for some $\mu > 0$
$$
\hat{h}_2(z, \hat{z}, x_0, x_1(x_0, z^*)) < h_2(z, x_0, x_1) 
(h_1(z, x_0) / \hat{h}_1(z, \hat{z}, x_0) + \mu).
$$
On the other hand, when $x_0 = (D \beta^{-1} \gamma^{-1} )^{1/2}$,
for some $\mu > 0$
$$
\hat{h}_1(z, \hat{z}, x_0) < h_1(z, x_0) (z / \hat{z} + \mu)
$$
and
$$
\hat{h}_2(z, \hat{z}, x_0, x_1(x_0, z^*)) > h_2(z, x_0, x_1) (\hat{z} / z - \mu).
$$
Since $\hat{h}_2$ is a continuous function, for large enough $D$ there's a solution
$(x_0, x_1)$ to (\ref{equ_C_V_equation}), such that 
$x_0 \in [\delta, (D \beta^{-1} \gamma^{-1} )^{1/2}]$. Due to
Lemma \ref{lem_jacobian}, for some $\mu > 0$ and any $Q \in \mathcal{B}(\mathbb{R}^2)$
$$
P_{(z, \hat{z})}((z_2, \hat{z}_2) \in Q \cap Z_V) \geq \mu \int_{Q \cap Z_V} d \lambda_V.
$$
To ensure the strict drift condition (\ref{equ_strict_drift}), it only
remains to show that
$$
\sup_{z \in \mathbb{R}^+, \hat{z} \in \mathbb{R}^+}
(E_{(z,\hat{z})} [ V_V(z_1, \hat{z}_1)] - V_V(z, \hat{z})) < \infty.
$$
It easily follows from (\ref{equ_RV_estimate_simple_z_zhat}) and
(\ref{equ_RV_estimate_any_z_log_diff_n}). Thus, the correctness of
(\ref{equ_strict_drift}) is established. Theorem \ref{thm_drift_criteria}
implies that $\Phi$ is indeed recurrent.

\textbf{Step 3.} Now, (\ref{equ_strict_recurrence}) and $\phi$-irreducibility
are checked in a similar
fashion to the Step 3 of Theorem \ref{thm_stat_distr_vgarch}.
Applying the Theorem \ref{thm_inv_measure_existence} and (\ref{equ_strict_recurrence})
yields the desired result.$\Box$


\subsection{Proof of the case i) of Theorems \ref{thm_lln_noninv_egarch} and
\ref{thm_lln_noninv_vgarch}}
The proofs for these laws of large numbers will also require some
general Markov chain technique. Essentially, they are corollaries of
theorems \ref{thm_stat_distr_egarch} and \ref{thm_stat_distr_vgarch} and
of corresponding Markov chains being positive Harris recurrent.
\begin{definition}
The set $A$ is called Harris recurrent if for all $x \in A$
$$
Q(x,A) = P(\Phi_t \in A \; \mbox{infinitely often} \, | \Phi_0 = x) = 1.
$$
The Markov chain $\{\Phi_t\}$ is called Harris recurrent if any set
in $\mathcal{B}^+$ is Harris recurrent.
\end{definition}

We will proof the case i) of both theorems simultaneously. Consider the chain $\{(z_t, \hat{z}_t)\}$
defined either by (\ref{equ_egarch_evol_zhat}) or (\ref{equ_vgarch_evol}).
For such a chain, positivity was proven in theorems \ref{thm_stat_distr_egarch} and
\ref{thm_stat_distr_vgarch}. Aperiodicity follows from the fact that
$$
\cup_{n \in \mathbb{N}} C_E(n, n) = \mathcal{X}_E, \quad
\cup_{n \in \mathbb{N}} C_V(n, n) = \mathcal{X}_V
$$
and petiteness of $C_E(A, B)$, $C_V(A, B)$.
In the Step 3 of theorems \ref{thm_stat_distr_egarch} and \ref{thm_stat_distr_vgarch}
we showed that for $A$, $B$ large enough and any $x \in C$ $L(x, C) = 1$,
where $C$ is either $C_E(A, B)$ or $C_V(A, B)$.
 By \cite[Proposition 9.1.1]{MeynTweedie2009}, this
implies Harris recurrence. By \cite[Theorem 13.0.1]{MeynTweedie2009},
for any $x \in \mathcal{X}$
\begin{equation}\label{equ_distr_convergence}
\sup_{A \in \mathcal{B}(X)} |P(\Phi_n \in A | \Phi_0 = x) - \pi(A)| \to 0
\end{equation}
as $n \to \infty$. By \cite[Theorem 17.0.1]{MeynTweedie2009}, positive
Harris recurrence also implies LLN of the form
\begin{equation}\label{equ_sample_path_convergence}
n^{-1} \sum_{t = 1}^n f(\overline{\Phi}_t) \stackrel{a.s.} \to
E f(\overline{\Phi}_0)
\end{equation}
for any $f$, $E |f(\overline{\Phi}_0)| < \infty$, where $\{\overline{\Phi}_t\}$ is
a stationary version of $\{\Phi_t\}$. It remains to check that the
statement we need is a combination of (\ref{equ_distr_convergence}) and
(\ref{equ_sample_path_convergence}). Indeed, denote
$$
\mathcal{X}_{LLN} = \left \{x \in \mathcal{X} :
P \left ( \left . n^{-1} \sum_{t = 1}^n f(\Phi_t) \to E f(\overline{\Phi}_0) \right |
\Phi_0 = x \right ) = 1 \right \}.
$$
By (\ref{equ_sample_path_convergence}) and Fubini theorem applied to
$\overline{\Phi}_0$ and $\varepsilon_0, \ldots \varepsilon_{n-1}$
$$
P(\overline{\Phi}_0 \in \mathcal{X}_{LLN}) = 1,
$$
which implies ergodicity for the process $\{\overline{\Phi}_t\}$. By (\ref{equ_distr_convergence}) and
Borel-Cantelli lemma, for any $x \in \mathcal{X}$
$$
P \left ( \exists n \geq 0, \Phi_n \in \mathcal{X}_{LLN} | \Phi_0 = x \right ) = 1.
$$
Therefore, for any $x \in \mathcal{X}$
$$
P \left ( \left .n^{-1} \sum_{t = 1}^n f(\Phi_t) \to E f(\overline{\Phi}_0) \right |
 \Phi_0 = x \right ) = 1
$$
and for any random variable $\xi$ which is measurable with respect to 
$\sigma\{\Phi_t, \; t \leq 0\}$
$$
P \left ( \left .n^{-1} \sum_{t = 1}^n f(\Phi_t) \to E f(\overline{\Phi}_0)
\right | \Phi_0 = \xi \right ) = 1.
$$
Setting $\xi = (\overline{\sigma}_0^2, s^2)$, where $s^2$ is any
starting point for $\hat{\sigma}^2_t$, concludes the proof. $\Box$


\subsection{Proof of Theorem \ref{thm_lln_noninv_egarch} ii)}

We will use notation from the proof of case i). It suffices to check that 
for the chain $\{(z_t, d_t)\}$, defined in (\ref{equ_egarch_evol}) and any $\mu > 0$
and $(z, d) \in \mathcal{X}_E$
$$
P_{(z, d)}(|d_t| > \mu) \to 0, \quad t \to \infty.
$$

We will use Theorem 8.4.3 from \cite{MeynTweedie2009} in order to 
establish transience of the chain $\{(z_t, d_t)\}$. Denote
$$
V_E^B(z, d) = 1 - \frac{1}{\max(-\log |d|, 1)}.
$$
For brevity, denote
$$
x \vee y = \max(x, y).
$$
Let us now check that for some $\mu > 0$ and any
$(z,d)$, such that $V_E^B(z, d) > 1 - \mu$,
\begin{equation}\label{equ_VE_estimate_transiency}
E_{(z,d)} [V_E^B(z_1, d_1)] > V_E^B(z, d).
\end{equation}
Indeed, for $|d| < \exp\{-1\}$ and any $z \in \mathbb{R}$ due to Fatou's lemma
$$
\mathop{\lim \inf}_{|d| \to 0} \left ( (\log |d|)^2 \left (E_{(z,d)} [V_E^B(z_1, d_1)] - V_E^B(z, d) \right ) \right ) = 
$$
$$
\mathop{\lim \inf}_{|d| \to 0} \left ( (\log |d|)^2 \left (E_{(z,d)} 
\frac{- 1 - (\log |d|)^{-1} ((-\log |\beta d + e_d \varepsilon_0^* |) \vee
 1 ) }{(-\log |\beta d + e_d \varepsilon_0^*|) \vee 1} \right ) \right ) =
$$
$$
\mathop{\lim \inf}_{|d| \to 0}  \left (- \log |d| E_{(z,d)} \frac{(-\log |\beta + d^{-1} e_d \varepsilon_0^* |) \vee
 (1 + \log |d|) }{(-\log |d| -\log |\beta + d^{-1} e_d \varepsilon_0^*|) \vee 1} \right ) =
$$
$$ 
\mathop{\lim \inf}_{|d| \to 0}  E_{(z,d)} \frac{(-\log |\beta + d^{-1} e_d \varepsilon_0^* |) \vee
 (1 + \log |d|) }{(1 + (\log |d|)^{-1} \log |\beta + d^{-1} e_d \varepsilon_0^*|) \vee (-(\log |d|)^{-1})}
\geq
$$
$$
\mathop{\lim \inf}_{|d| \to 0}  E \left [ -\log |\beta + \varepsilon^*_0 d^{-1} e_d| \right ] = - \lambda^y_E > 0.
$$
Therefore, (\ref{equ_VE_estimate_transiency}) holds true 
for all small $|d|$ as needed.

As was shown during the proof of Step 2 of Theorem
\ref{thm_stat_distr_egarch}, the set $C_E(A, B)$ is petite for any 
$0 < A < \infty$, $0 < B < \infty$. During Step 3 of the same Theorem,
it was established that the chain $\{(z_k, d_k)\}$ is 
$\psi$-irreducible. The exact same proof of irreducibility works in the case $\lambda^E < 0$
too, except for the equality $L(x, C_E(A, B)) = 1$ which is not true anymore.
Instead, we'll verify directly that 
for large enough $A$ and $B$ and any $(z, d) \in \mathcal{X}_E$ there exists
$k$, such that $P_{(z,d)}((z_k, d_k) \in C_E(A,B)) > 0$. It suffices 
to check that for any $B > 1$ and large enough $A$
there exist $t, x_0, \ldots, x_{t-1}$, such that
$$
(h_t(z_0, x_0, \ldots, x_{t-1}), 
\hat{h}_t(d_0, x_0, \ldots, x_{t-1})) \in C_E(A, B).
$$
Indeed, for any $d \in \mathbb{R}$ there exist $x_0(d), x_1(d)$, such that
$$
\hat{h}_2(d, x_0(d), x_1(d)) = 1.
$$
Then for any $z, d$ 
$$
\hat{h}_{2 k}(d, x_0(d), x_1(d), x_0(1), x_1(1), \ldots, x_0(1), x_1(1)) = 1
$$
and 
$$
h_{2 k}(z, x_0(d), x_1(d), x_0(1), x_1(1), \ldots, x_0(1), x_1(1)) \to
\beta^{-2} (\alpha \beta + x_0(1) \beta + x_1(1)), \quad t \to \infty.
$$
It remains to choose any $A > \beta^{-2} (\alpha \beta + x_0(1) \beta + x_1(1))$.
Hence, the chain $\{(z_k, d_k)\}$ is $\psi$-irreducible.
It's easy to check now that $C_E(A, B)$ and $\mathcal{X}_E \backslash C_E(A, B)$ 
are both in $B^+(\mathcal{X}_E)$.

Due to (\ref{equ_VE_estimate_transiency})
and Theorem 8.4.2 from \cite{MeynTweedie2009}, the chain $\{(z_t, d_t)\}$
is transient. According to the Theorem  8.3.5 from \cite{MeynTweedie2009}
and the fact that $C_E(A, B)$ is petite, it is also uniformly transient.
Finally, since for arbitrary $A > 0$
$$
P_{(z, d)}(|d_t| > \mu) \leq P ((z_t, d_t) \in C_E(A, \mu^{-1})) + 
P_{(z, d)} (|z_t| > A) \to P (|\overline{z}_t| > A), \quad t \to \infty,
$$
we get
$$
P_{(z, d)}(|d_t| > \mu) \to 0, \quad t \to \infty. \Box
$$


\subsection{Proof of Theorem \ref{thm_lln_noninv_vgarch} ii)}

We need to check that 
for the chain $\{(z_t, d_t)\}$, defined in (\ref{equ_vgarch_evol}) and any $\mu > 0$
and $(z, \hat{z}) \in \mathcal{X}_V$
$$
P_{(z, \hat{z})}(|\hat{z}_t - z_t| > \mu) \to 0, \quad t \to \infty.
$$
Fix $m \in \mathbb{R}^+$, $n \in \mathbb{N}$ and denote
$$
V_V^B(z, \hat{z}) = 1 - \frac{1}{(m R_n(z, \hat{z}) - \sqrt{z}) \vee 1},
$$
with $R_n(z, \hat{z})$ as defined during the Step 1 of the proof of Theorem
\ref{thm_stat_distr_vgarch}. 
Let us now check that for some $\mu > 0$, $n \in \mathbb{N}$, 
$m \in \mathbb{R}^+$ and any
$(z, \hat{z})$, such that $V_V^B(z, \hat{z}) > 1 - \mu$,
\begin{equation}\label{equ_VV_estimate_transiency}
E_{(z,d)} [V_V^B(z_1, \hat{z}_1)] > V_V^B(z, \hat{z}).
\end{equation}

First, recall that due to (\ref{equ_RV_estimate_any_z_log_diff_n}) 
for any $(z, d) \in \mathcal{X}_V$,
$n \in \mathbb{N}$
\begin{equation}\label{equ_estimate_Rn_Rn1}
|R_n(z_1, \hat{z}_1) + \log|\hat{z}_1 - z_1| - 
R_n(z, \hat{z}) - \log|\hat{z} - z|| \leq 2 M^* n,
\end{equation}
with $M^*$ as defined by (\ref{equ_M_star}).

It's easy to see that for $x \geq 1$ and $y \in \mathbb{R}$
\begin{equation}\label{equ_one_over_R_n_estimate}
\frac{1}{(x+y) \vee 1} - \frac{1}{x} + \frac{y}{x^2} \leq
\frac{y^2}{x^2 ((x+y) \vee 1)} \leq \frac{y^2}{x^2}.
\end{equation}
Put 
$$
d(z) = (\sqrt{z} - A) \vee 0,
$$
we will select $A$ in what follows.
For $(z, \hat{z})$ such that $V_V^B(z, \hat{z}) > 0$ we have
$$
V_V^B(z_1, \hat{z}_1) \geq 1 - \frac{1}{(m R_n(z_1, \hat{z}_1) - d(z) \vee \sqrt{z_1}) \vee 1} = 
V_V^B(z, \hat{z}) + S_n(z, \hat{z}, z_1, \hat{z}_1) + T_n(z, \hat{z}, z_1, \hat{z}_1),
$$
where
$$
T_n(z, \hat{z}, z_1, \hat{z}_1) = \frac{(m R_n(z_1, \hat{z}_1) - d(z) \vee \sqrt{z_1}) \vee 1 - 
m R_n(z, \hat{z}) + \sqrt{z}}{(m R_n(z, \hat{z}) - \sqrt{z})^2},
$$
$$
S_n(z, \hat{z}, z_1, \hat{z}_1) = \frac{1}{m R_n(z, \hat{z}) - \sqrt{z}} -
\frac{1}{(m R_n(z_1, \hat{z}_1) - d(z) \vee \sqrt{z_1}) \vee 1} - T_n(z, \hat{z}, z_1, \hat{z}_1).
$$
Let's estimate $S_n$ and $T_n$ separately. First, using  
(\ref{equ_one_over_R_n_estimate}), (\ref{equ_estimate_Rn_Rn1}) and
the inequality
\begin{equation}\label{equ_sqrt_estimate}
\sqrt{y} - \sqrt{x} \leq \sqrt{y - x}
\end{equation}
which is valid for any $y \geq x > 0$, we get 
$$
(m R_n(z, \hat{z}) - \sqrt{z})^2 S_n(z, \hat{z}, z_1, \hat{z}_1) \leq 
\frac{(m R_n(z_1, \hat{z}_1) - d(z) \vee \sqrt{z_1} - 
m R_n(z, \hat{z}) + \sqrt{z})^2}{(m R_n(z_1, \hat{z}_1) - d(z) \vee \sqrt{z_1}) \vee 1} \leq
$$
$$
\frac{2 m^2 (R_n(z_1, \hat{z}_1) - R_n(z, \hat{z}))^2 + 2(\sqrt{z} - d(z) \vee \sqrt{z_1})^2}
{(m R_n(z_1, \hat{z}_1) - d(z) \vee \sqrt{z_1}) \vee 1} \leq
$$
\begin{equation}\label{equ_sn}
\frac{4 m^2 (4 (M^*)^2 n^2 + (\log|\hat{z} - z| - \log|\hat{z}_1 - z_1|)^2) + 
2 A^2 + 2 \gamma (\varepsilon_0 - \delta)^2}
{(m R_n(z_1, \hat{z}_1) - d(z) \vee \sqrt{z_1}) \vee 1}.
\end{equation}
According to (\ref{equ_RV_estimate_any_z_log_diff_interm}) and Lemma 
\ref{lem_log_estimate} (\ref{equ_log_estimate_f}), the nominator
(\ref{equ_sn}) is uniformly integrable for $(z, \hat{z}) \in \mathcal{X}_V$,
and its denominator tends to infinity for a fixed $\epsilon_0$ as 
$$
m R_n(z, \hat{z}) - \sqrt{z} \to +\infty.
$$
Due to dominated convergence theorem,
$$
E_{(z, \hat{z})} \Bigl [(m R_n(z, \hat{z}) - \sqrt{z})^2 S_n(z, \hat{z}, z_1, \hat{z}_1) \Bigr ] \to 0
$$
as $m R_n(z, \hat{z}) - \sqrt{z} \to +\infty$.

For $T_n$ we will consider 2 cases.

1) $z \leq e^{\sqrt{n}}$. Similarly to the estimation of
the right-hand side of (\ref{equ_RV_estimate_z_leq_exp_sqrt_n}), when $n$
is large enough, $\mu > 0$ is small enough, for any $(z, \hat{z})$ such that
$$
R_n(z, \hat{z}) > \mu^{-1},
$$
we have
$$
E_{(z, \hat{z})} [R_n(z_1, \hat{z}_1) - R_n(z, \hat{z})] \geq \lambda_V^y / 2.
$$
Therefore, for any such $n$, $\mu$
$$
m > (\lambda_V^y / 2)^{-1} (\sqrt{\gamma} E |\varepsilon_0 - \delta| + A),
$$
we have
$$
E_{(z, \hat{z})} T_n(z, \hat{z}, z_1, \hat{z}_1) \geq \frac{m E [R_n(z_1, \hat{z}_1) - R_n(z, \hat{z})] -
E_{(z, \hat{z})} [d(z) \vee \sqrt{z_1}]  + \sqrt{z}}{(m R_n(z, \hat{z}) - \sqrt{z})^2}  \geq
$$
$$
\frac{m E_{(z, \hat{z})} [R_n(z_1, \hat{z}_1) - R_n(z, \hat{z})] -
d(z) \vee (\sqrt{\beta z} + \sqrt{\gamma} E |\varepsilon_0 - \delta|) +
 \sqrt{z}}{(m R_n(z, \hat{z}) - \sqrt{z})^2} \geq
$$
\begin{equation}\label{equ_tn_part1}
\frac{m \lambda_V^y / 2 - \sqrt{\gamma} E |\varepsilon_0 - \delta|}
{(m R_n(z, \hat{z}) - \sqrt{z})^2} > 0.
\end{equation}

2) $z > e^{\sqrt{n}}$. According to (\ref{equ_RV_estimate_any_z_log_diff}), 
(\ref{equ_estimate_Rn_Rn1}), Lemma \ref{lem_log_estimate} (\ref{equ_log_estimate_f})
 and (\ref{equ_sqrt_estimate})
$$
E T_n(z, \hat{z}, z_1, \hat{z}_1) \geq E \frac{m R_n(z_1, \hat{z}_1) - d(z) \vee \sqrt{z_1}  - 
m R_n(z, \hat{z}) + \sqrt{z}}{(m R_n(z, \hat{z}) - \sqrt{z})^2} \geq
E \frac{- m M - d(z) \vee (\sqrt{\beta z} + \sqrt{\gamma} |\varepsilon_0 - \delta|) +
 \sqrt{z}}{(m R_n(z, \hat{z}) - \sqrt{z})^2} \geq
$$
\begin{equation}\label{equ_tn_part2}
\frac{- m (2 M^* + M) - E (\sqrt{\gamma} |\varepsilon_0 - \delta|)}{(m R_n(z, \hat{z}) - \sqrt{z})^2}.
\end{equation}
In view of (\ref{equ_sn}), (\ref{equ_tn_part1}) and (\ref{equ_tn_part2}),
(\ref{equ_VV_estimate_transiency}) is established.

XXX Add the proof of $psi$-irreducibility.



\subsection{Auxilliary results}
The following lemmas are technical tools used in proofs of Theorems
\ref{thm_stat_distr_egarch} and \ref{thm_stat_distr_vgarch}.
\begin{lemma}\label{lem_log_estimate}
Assume $\xi$ is a random variable with a distribution absolutely continuous 
with respect to Lebesgue measure with density $h(x)$, such that
$$
\sup_{x \in \mathbb{R}} (|x| + 1) h(x) < \infty.
$$
Then
\begin{equation*}\tag{a}\label{equ_log_estimate_b}
E (-\log|1 + q \xi|)^+ \to 0, \quad q \to 0.
\end{equation*}
\begin{equation*}\tag{b}\label{equ_log_estimate_c}
E (-\log|1 + q \xi|)^+ \to 0, \quad q \to \infty.
\end{equation*}
and for any $C < \infty$ and $\rho \geq 1$ 
\begin{equation*}\tag{c}\label{equ_log_estimate_f}
\sup_{|q| \leq C, 0 < r \leq C} E |\log|1 + q \xi - r \xi^2| |^{\rho} < \infty.
\end{equation*}
\end{lemma}


\textbf{Proof.}

(a) If $q = 0$ the left-hand side is $0$. If $q > 0$, due to Condition
\ref{cond_regul} for any $\mu \in (0, 1)$
$$
E (-\log|1 + q \xi|)^+ = - \int_{-2/q}^0 \log(|1 + q x|) h(x) d x
\leq
$$
$$
-\log(1 - \mu) - \int_{-2 + \mu}^{-\mu} \log(|1 + y|) h(y/q) q^{-1} d y \leq 
-\log(1 - \mu) - M \int_{-2 + \mu}^{-\mu} \log(1 + y) h_q(y) d y,
$$
where $0 \leq h_q(y) \leq |y|^{-1}$ for all $y \in \mathbb{R} \backslash \{0\}$.
Noting that $h_q(y)$ is uniformly bounded on $[-2 + \mu, -\mu]$,
$\log(1 + y)$ is integrable over the same interval and that
$$
\int_{-2 + \mu}^{-\mu} h_q(y) d y \to 0, \quad q \to 0
$$
concludes the proof (using a slightly modified dominated convergence
theorem). The case $q < 0$ is similar.

(b) For $q > 0$ and $\mu \in (0, 1)$
$$
E (-\log|1 + q \xi|)^+ =
- \log(1 - \mu) - \int_{-2 + \mu}^{-\mu} \log(|1 + y|) h(y/q) q^{-1} d y
\leq
$$
$$
-\log(1 - \mu) - \left [ \sup_{x \in \mathbb{R}} h(x) \right ]
q^{-1} \int_{-2 + \mu}^{-\mu} \log(|1 + y|) d y \to
-\log(1 - \mu), \quad q \to \infty.
$$
Since $\mu$ may be arbitrarily small, the proof for $q > 0$ is thus complete.
The case $q < 0$ is similar.

(c) For any $|q| \leq C, 0 < r \leq C$ the polynom $1 + q x - r x^2$ has
2 roots, denote them $a_1^{-1}$ and $a_2^{-1}$. Then 
$$
|1 + q \xi - r \xi^2| = |1 - a_1 \xi| 
|1 - a_2 \xi|,
$$
and it's easy to check that 
$$
\sup_{|q| \leq C, 0 < r \leq C} |a_i|  < \infty, \quad i = 1, 2.
$$
It suffices to show that
$$
\sup_{|q| \leq C, 0 < r \leq C} E |\log |1 + q \xi - r \xi^2| |^{\rho} < \infty.
$$
Since 
$$
E |\log |1 + q \xi - r \xi^2| |^{\rho} \leq 2^{\rho} 
E |\log |1 - a_1 \xi^2||^{\rho} + 2|^{\rho} E |\log |1 - a_2 \xi^2||^{\rho},
$$
it is enough to check that
$$
\sup_{|a| \leq C} E |\log|1 - a \xi | |^{\rho} < \infty.
$$
Indeed, for $a>0$ due to Condition \ref{cond_regul}
$$
E |\log|1 - a \xi ||^{\rho} = 
\int_{|1 - a x| \leq 2} |\log |1 - a x||^{\rho} h(x) d x +
 \int_{|1 - a x| > 2} |\log |1 - a x||^{\rho} h(x) d x \leq
$$
$$
\int_{-2}^2 |\log |y||^{\rho} a^{-1} h((1 -y) a^{-1}) d y + 
\sup_{x \geq 2} (x^{-1} |\log x|^{\rho} ) \int (1 + a |x|) h(x) d x \leq 
$$
$$
M \int_{-2}^2 |\log |y||^{\rho} |1 -y|^{-1} d y + 
\sup_{x \geq 2} (x^{-1} |\log x|^{\rho} ) (1 + a E |\xi|) < \infty,
$$
where, once again, $M = \sup_{x \in \mathbb{R}} \left ( (|x| + 1) h(x) \right )$. 
The case $a < 0$ is similar, and $a = 0$ is obvious. $\Box$

\begin{lemma}\label{lem_jacobian}
Let $\mathcal{T}$ be an arbitrary set and
$f_{\tau} : \mathbb{R}^k \to \mathbb{R}^k$ be a
continuously differentiable function for any $\tau \in \mathcal{T}$.
Assume also that for any $\mathbf{z} \in Z \subseteq \mathbb{R}^k$ the equation
$$
f_{\tau}(x_1, \ldots, x_k) = \mathbf{z},
$$
has a solution
$$
\mathbf{x}(\tau, \mathbf{z}) = (x_1(\tau, \mathbf{z}), \ldots,
x_k(\tau, \mathbf{z})),
$$
which is uniformly bounded for $\tau \in \mathcal{T}$, $\mathbf{z} \in Z$.

Also assume $\xi_1, \ldots, \xi_k$ to be i.i.d. random
variables with and absolutely continuous distribution whose density
is uniformly separated from $0$ on 
$$
\{x_j(\tau, \mathbf{z}), \; j = 1, \ldots, k, \tau \in \mathcal{T}, \mathbf{z} \in Z\}.
$$
Put $\boldsymbol{\xi} = (\xi_1, \ldots, \xi_k)$.
Then for some $\mu > 0$ and any $C \subseteq Z$, $\tau \in \mathcal{T}$
$$
P \left (f_{\tau} (\boldsymbol{\xi}) \in C\right ) \geq
\mu \lambda(C),
$$
where $\lambda$ is standard Lebesgue measure on $\mathbb{R}^k$.
\end{lemma}
\textbf{Proof.}
Denote $u_{\delta}(\mathbf{x})$ a $L_2$-ball in $\mathbb{R}^k$
with a center in $\mathbf{s}$ and radius $\delta > 0$, then
for small enough $\nu > 0$ and any $\tau \in \mathcal{T}$, $\mathbf{z} \in Z$
\begin{equation}\label{equ_density_estimate}
P \left (f_{\tau} (\boldsymbol{\xi}) \in u_{\delta}(\mathbf{z}) \right ) \geq
P (\boldsymbol{\xi} \in u_{\nu \delta}(\mathbf{x} (\tau, \mathbf{z}))) \geq
\nu \lambda(u_{\nu \delta}(\mathbf{x} (\tau, \mathbf{z})))) =
\nu^{k + 1} \lambda(u_{\delta}(\mathbf{z})).
\end{equation}
The inequality (\ref{equ_density_estimate}) implies that for some $\mu > 0$
and any parallelepiped $C \subset \mathbb{R}^k$
$$
P \left (f_{\tau} (\boldsymbol{\xi}) \in C \right ) \geq \mu \lambda(C).
$$
The application of monotone class theorem concludes the proof.$\Box$

The following is a direct corollary of Lemma \ref{lem_jacobian}:
\begin{lemma}\label{lem_jacobian_more}
In the conditions of Lemma \ref{lem_jacobian}, assume that
$$
\mathcal{T} \subseteq \mathcal{T}_1 \times \mathcal{T}_2,
$$
where $\mathcal{T}_1$ is an arbitrary set and $\mathcal{T}_2$ is a
measurable space and for any $\tau_1 \in \mathcal{T}_1$
$$
\{\tau_2 \in \mathcal{T}_2: (\tau_1, \tau_2) \in \mathcal{T}\}
$$
is measurable.
Also assume that for some random element $\chi$ in $\mathcal{T}_2$,
independent with $\boldsymbol{\xi}$,
$$
P((\tau_1, \chi) \in \mathcal{T}) \geq \rho > 0
$$
for any $\tau_1 \in \mathcal{T}_1$.
Then for some $\mu > 0$ and any $C \subseteq Z$, $\tau_1 \in \mathcal{T}_1$
$$
P \left (f_{(\tau_1, \chi)} (\boldsymbol{\xi}) \in C\right ) \geq
 \mu \lambda(C).
$$
\end{lemma}


\section{Acknowledgements}

The author is grateful for the support of AHL Research during work on this paper.
I would also like to thank
Dr. Jeremy Large, Professor Neil Shephard, Dr. Kevin Sheppard
(all of Oxford-Man Institute of Quantitative Finance), Professor
Anders Rahbek (Copenhagen) and Professors Michael Boldin, Yuri Tyurin and
Valeri Tutubalin (all of Moscow State University) for helpful discussions.

The author also thanks anonymous referees for their suggestions which helped
to improve the presentation.


\pagebreak

\bibliographystyle{plainnat}

\end{document}